\documentclass{amsart}

\usepackage{amsmath,amssymb,amsthm}
\usepackage{xcolor}
\usepackage{aliascnt}
\usepackage{cleveref}
\usepackage[shortlabels]{enumitem}
\usepackage{needspace}
\numberwithin{equation}{section}

\newtheorem{theorem}{Theorem}[section]

\newaliascnt{proposition}{theorem}
\newtheorem{proposition}[proposition]{Proposition}
\aliascntresetthe{proposition}
\crefname{proposition}{proposition}{propositions}
\Crefname{proposition}{Proposition}{Propositions}

\newaliascnt{lemma}{theorem}
\newtheorem{lemma}[lemma]{Lemma}
\aliascntresetthe{lemma}
\crefname{lemma}{lemma}{lemmas}
\Crefname{lemma}{Lemma}{Lemmas}

\newaliascnt{definition}{theorem}
\newtheorem{definition}[definition]{Definition}
\aliascntresetthe{definition}
\crefname{definition}{definition}{definitions}
\Crefname{definition}{Definition}{Definitions}

\newaliascnt{remark}{theorem}
\newtheorem{remark}[remark]{Remark}
\aliascntresetthe{remark}
\crefname{remark}{remark}{remarks}
\Crefname{remark}{Remark}{Remarks}

\newaliascnt{example}{theorem}
\newtheorem{example}[example]{Example}
\aliascntresetthe{example}
\crefname{example}{example}{examples}
\Crefname{example}{Example}{Examples}

\newaliascnt{corollary}{theorem}
\newtheorem{corollary}[corollary]{Corollary}
\aliascntresetthe{corollary}
\crefname{corollary}{corollary}{corollaries}
\Crefname{corollary}{Corollary}{Corollaries}

\newcommand{\Cst}{\mathrm{C}^*}
\newcommand{\bN}{\mathbb{N}}

\newcommand{\bC}{\mathbb{C}}

\usepackage{biblatex}
\begin{document}

\title[Duality of Hopf $\Cst$-algebra actions]{Duality of finite-dimensional Hopf $\Cst$-algebra actions via inclusions}
\begin{abstract}
We study actions of finite-dimensional Hopf $\Cst$-algebras on a unital
$\Cst$-algebra $A$. We first define the Rokhlin property and approximate
representability directly for an action $H\curvearrowright A$ and
characterize these properties, under appropriate hypotheses, in terms
of the corresponding properties of the canonical conditional expectation
$E\colon A\to A^H$. Based on this characterization, we establish the
following duality for a saturated action on a simple unital
$\Cst$-algebra $A$ such that $A\rtimes H$ is simple:
\begin{enumerate}
\item $H\curvearrowright A$ has the Rokhlin property if and only if
the dual action $H^{\circ}\curvearrowright A\rtimes H$ is approximately
representable.
\item $H\curvearrowright A$ is approximately representable if and only
if the dual action $H^{\circ}\curvearrowright A\rtimes H$ has the Rokhlin
property.
\end{enumerate}
For separable $A$, we also identify these properties with the
corresponding spatial formulations introduced by Barlak--Szab{\'o}--Voigt,
assuming saturation in the Rokhlin comparison. We present examples of
Hopf $\Cst$-algebra actions both with and without the Rokhlin property.
\end{abstract}
\author{Anshu}
\address{Department of Mathematics, University of Ulsan,
Ulsan 44610, Republic of Korea}
\email{anshu@ulsan.ac.kr}

\author{Hyun Ho Lee}
\address{Department of Mathematics, University of Ulsan,
Ulsan 44610, Republic of Korea}
\email{hadamard@ulsan.ac.kr}
\date{}
\subjclass[2020]{Primary 46L55; Secondary 46L05, 46L37, 16T05}
\keywords{Rokhlin property, approximate representability, Hopf
$\Cst$-algebra action, conditional expectation, Watatani index,
basic construction, duality}
\thanks{This work was supported by~the National Research Foundation of Korea
(NRF) grant funded by the Korean government (MSIT)
(NRF-RS-2025-00557662).}

\maketitle

\section{Introduction}

M.~Izumi studied the Rokhlin property for finite group actions on
$\Cst$-algebras in \cite{Izumi2004I,Izumi2004II}. For actions of finite
abelian groups on separable unital $\Cst$-algebras, he proved that the Rokhlin property and approximate
representability are dual to each other
\cite[Lemma~3.8]{Izumi2004I}. Osaka, Kodaka, and Teruya extended these
notions to inclusions of index-finite type and introduced a Rokhlin property
for finite-dimensional Hopf $\Cst$-algebra actions
\cite[Definition~4.4]{OKT-2009}. Kodaka and Teruya subsequently showed
that this earlier formulation forces the dual Hopf $\Cst$-algebra to be
commutative \cite[Corollary~11.5, p.~369]{KT-2015}.

For finite-dimensional Hopf $\Cst$-algebra coactions, Kodaka and Teruya
introduced the Rokhlin property through approximate representability
of the dual coaction \cite[Definition~5.1]{KT-2015}. Importantly, they
also proved an intrinsic characterization: a coaction has the Rokhlin
property precisely when its central sequence algebra contains a
projection whose Haar average is $1/\dim H$. These results were established in \cite[Theorem~5.12]{KT-2015} under the assumption of saturation and in
\cite[Corollary~6.4]{KT-2015} without this assumption.

From a different perspective, Barlak, Szab{\'o}, and Voigt introduced
the spatial Rokhlin property for coexact compact quantum group actions
and spatial approximate representability for discrete quantum group
actions, using equivariantly sequentially split maps and braided
tensor products.
They defined the two properties independently and proved their duality
\cite[Definitions~4.1 and~4.7, and Theorem~4.12]{BSV2017}; for finite
quantum groups, they recover the Kodaka--Teruya notions
\cite[Remark~6.2]{BSV2017}. Their framework is particularly effective
for proving structural results and developing general theory, but
verifying the required equivariant splitting for a specific action
can be difficult.

In contrast to Kodaka and Teruya's duality-based definition of the
Rokhlin property, we introduce direct formulations of the Rokhlin
property and approximate
representability for finite-dimensional Hopf $\Cst$-algebra actions
and characterize these properties through the fixed-point inclusion
$A^H\subset A$. More precisely, under appropriate hypotheses, each
property of $H\curvearrowright A$ is equivalent to the corresponding
property of the canonical conditional expectation $E\colon A\to A^H$.
Theorem~\ref{thm:main} shows that a saturated action is
approximately representable if and only if $E$ is approximately
representable, provided $A^H$ is simple; the implication from $E$ to
the action needs no simplicity assumption. This extends
\cite[Proposition~3.17]{LO2019}, which treats finite abelian groups
with $A$ and $A\rtimes_\alpha G$ simple, to Hopf $\Cst$-algebras and
in particular to nonabelian groups. Theorem~\ref{thm:Hopf-incl} gives
the corresponding statement for the Rokhlin property when $A$ is
simple. Our Rokhlin formulation requires centrality only of the image
of the distinguished projection and hence avoids the commutativity
restriction of the earlier Osaka--Kodaka--Teruya definition.

The duality is therefore determined entirely by the inclusion. The key identity is
$F=\widehat E$: the fixed-point expectation for the dual action on
$A\rtimes H$ is the dual conditional expectation of $E$ arising from
the Jones basic construction
(Proposition~\ref{prop:crossed-product-basic-construction}).
Combining the two characterizations with the Lee--Osaka duality for
conditional expectations \cite[Theorem~3.16]{LO2019} proves
Theorem~\ref{Rok-AR-duality}: the Rokhlin property and approximate
representability are exchanged by passage to the dual action. This
proof does not invoke the existing dynamical duality theorems for
coactions or spatial quantum group actions. It explains this exchange
through the Jones basic construction.

We then compare our notions with the spatial ones.
Theorem~\ref{thm:RP-SRP} shows that, for saturated actions on separable
$A$, our Rokhlin property agrees with the spatial Rokhlin property, and
Theorem~\ref{thm:AR-SAR} shows the same for approximate
representability, without saturation or simplicity. Consequently, the
duality theorem of Barlak, Szab{\'o}, and Voigt
\cite[Theorem~4.12]{BSV2017} also applies to our notions and gives an
alternative route to Theorem~\ref{Rok-AR-duality} in the separable
case (see the remark following its proof). We apply the results to the
non-group type crossedproduct inclusion of
\cite[Proposition~5.8 and Example~5.9]{LO2026WT} and to a
Kac--Paljutkin UHF action, and observe some permanence consequences.
Finally, we use the projection in the definition of approximate
representability of $E$ to obtain an obstruction
(Theorem~\ref{thm:obstruction}), and combine it with duality in
Example~\ref{ex:jiang-su}.

This article is organized as follows. In Section~\ref{sec:prelim}, we collect
the terminology needed later and explain how an action of a Hopf
$\Cst$-algebra on a $\Cst$-algebra $A$ corresponds to a coaction of its
dual on $A$. The main results are proved in Section~\ref{sec:main}, and examples and applications are discussed in Section~\ref{sec:applications}.
Throughout, unless stated otherwise, $A$ denotes a unital
$\Cst$-algebra, $G$ denotes a classical group, and all Hopf
$\Cst$-algebras are finite-dimensional.


\section{Preliminaries}\label{sec:prelim}

\subsection{Hopf algebras and actions}
We use the following conventions; see \cite[Propositions~2.1 and~2.3]{SP1994}.

\begin{definition}
    A finite-dimensional \emph{Hopf $\Cst$-algebra} is a tuple $(H, \Delta_H, S_H, \varepsilon_H)$, where $H$ is a finite-dimensional unital $\Cst$-algebra, $\Delta_H$ and $\varepsilon_H$ are unital $\Cst$-algebra homomorphisms, and $S_H$ is a linear map satisfying the following conditions.
    \begin{enumerate}
        \item (comultiplication) $\Delta_H: H \to H \otimes H$ is coassociative, that is, 
        \[
        (\Delta_H \otimes \operatorname{id})\Delta_H = (\operatorname{id} \otimes \Delta_H) \Delta_H ; \]
        \item (counit) $\varepsilon_H: H \to \bC$ such that 
        \[
        (\varepsilon_H \otimes \operatorname{id}) \Delta_H = \operatorname{id} = (\operatorname{id} \otimes \varepsilon_H)\Delta_H;
        \]
        \item (antipode) The linear map $S_H:H\to H$ satisfies, for every $h\in H$,
        \[
        m_H(S_H\otimes\operatorname{id}_H)\Delta_H(h)
        =\varepsilon_H(h)1_H
        =m_H(\operatorname{id}_H\otimes S_H)\Delta_H(h)
        \]
        
    \end{enumerate}
where $m_H$ denotes multiplication in $H$.
\end{definition}

In the finite-dimensional $\Cst$-setting, the antipode is automatically
unital and antimultiplicative, with $S_H^2=\operatorname{id}_H$ and
$S_H(h^*)=S_H(h)^*$ for $h\in H$; see
\cite[Proposition~2.1(2)]{SP1994}. Thus $H$ is of Kac type, and these
properties are consequences rather than additional hypotheses.

For a Hopf algebra, there is an associated dual Hopf $\Cst$-algebra $(H^\circ, \Delta^\circ, S^\circ, \varepsilon^\circ)$; see \eqref{dual-H}. We recall some basic terminology for Hopf algebras and will use the following properties of their structure maps, whose proofs can be found in standard texts on quantum groups (see \cite{timmermann, SP1994, VD1997}).

 \begin{lemma}\label{prop: Hopf}
     Let $(H, \Delta, S, \varepsilon)$ be a finite-dimensional Hopf $\Cst$-algebra. Then there is a distinguished projection $e_H \in H$ such that 
         $$h e_H = e_H h = \varepsilon (h)e_H,  \ \varepsilon(e_H)=1,  \ S(e_H) = e_H.$$
 \end{lemma}

The action of a Hopf $\Cst$-algebra on a unital $\Cst$-algebra $A$ is defined as follows \cite[Definition~2.4]{SP1994}.

\begin{definition}\cite[Definition~2.4]{SP1994}\label{def:H-action}
 A bilinear map $\cdot : H \times A \to A$ is a \emph{left action} of $H$ on $A$ if the following conditions hold for all $h,k \in H$, $a,b \in A$.
\begin{align}
1_{H} \cdot a & = a, \tag{i} \\
h \cdot 1_A &= \varepsilon(h) 1_A, \tag{ii} \\
hk \cdot a &= h \cdot(k\cdot a), \tag{iii} \\
h \cdot ab &= \sum (h_{(1)} \cdot a )(h_{(2)} \cdot b), \tag{iv} \\
(h\cdot a)^* &= S(h^*)\cdot a^* \tag{v}
\end{align}
where $\Delta(h) = \sum h_{(1)}\otimes h_{(2)}$ is the Sweedler-type notation for comultiplication.
\end{definition}

 For the general framework of coactions of Hopf $\Cst$-algebras, we refer to Baaj and Skandalis \cite[Definition~0.2]{BS1993}. In the finite-dimensional unital setting considered here, we use the following definition of Kodaka and Teruya \cite[Definition~2.4]{KT2009}. Below, we recall the correspondence between an action of $H$ on $A$ and a coaction of its dual $H^\circ$ on $A$.

\begin{definition}\cite[Definition~2.4]{KT2009}\label{def:coaction}
A \emph{coaction} of a Hopf $\Cst$-algebra $(H,\Delta)$ on $A$ is a unital injective $*$-homomorphism $\rho\colon A  \to   A\otimes H$ such that
\begin{equation*}
(\rho \otimes \operatorname{id})\circ \rho=(\operatorname{id} \otimes \Delta) \circ \rho,
\end{equation*}
together with the counit condition $(\operatorname{id}\otimes\varepsilon) \circ \rho = \operatorname{id}_A$.
\end{definition}

For the remainder of the article, let us fix a basis $\{h_i\}_{i=1}^n$ of the finite-dimensional $\Cst$-algebra $H$ and the dual basis 
$\{h^i\}_{i=1}^n$ of $H^\circ$, with the dual pairing $\langle h_i, h^j\rangle  = h^j(h_i)=\delta_{ij}$. The dual Hopf $\Cst$-algebra structure on $H^\circ$ is given by \cite[Proposition~2.3]{SP1994}
\begin{equation}\label{dual-H}
\begin{aligned}
 (\phi\psi)(a)&=(\phi\otimes\psi)\Delta_H(a), &
 (\Delta^{\circ}\phi)(a\otimes b)&=\phi(ab),\\
 \varepsilon^{\circ}(\phi)&=\phi(1_H), &
 1_{H^\circ}&=\varepsilon_H,\\
 S^{\circ}(\phi)(a)&=\phi(S_H(a)), &
 \phi^*(a)&=\overline{\phi(S_H(a^*))}.
\end{aligned}
\end{equation}
where $a,b\in H$ and $\phi,\psi\in H^\circ$. The following correspondence
uses the conventions of \cite[Definitions~2.2 and~2.4]{KT2009}.
We include the verification to make the coproduct and involution
conventions explicit. Now we will write $K=H^\circ$ for convenience.

\begin{proposition}\label{thm: act-coact}
 Left actions of $H$ on $A$ in
Definition~\ref{def:H-action} correspond bijectively to coactions
$\rho\colon A\to A\otimes K$ in Definition~\ref{def:coaction} and the following assignments are mutually inverse to each other
\[
 \rho(a)=\sum_i(h_i\cdot a)\otimes h^i,
 \qquad
 h\cdot a=(\operatorname{id}_A\otimes\operatorname{ev}_h)\rho(a),
 \quad \operatorname{ev}_h(\phi)=\phi(h).
\]
\end{proposition}

\begin{proof}
For $h\in H$, write $F_h=\operatorname{id}_A\otimes
\operatorname{ev}_h$. These slices separate points of $A\otimes K$;
likewise, the slices $\operatorname{id}_A\otimes
\operatorname{ev}_h\otimes\operatorname{ev}_k$ separate points of
$A\otimes K\otimes K$.

Suppose first that $H$ acts on $A$, and define
$\rho(a)=\sum_i(h_i\cdot a)\otimes h^i$. The dual-basis identity
gives $F_h\rho(a)=h\cdot a$. Since $1_K=\varepsilon_H$ and
$\varepsilon_K(\phi)=\phi(1_H)$, the unit rules give
\[
\begin{aligned}
\rho(1_A)&=\sum_i\varepsilon_H(h_i)1_A\otimes h^i
           =1_A\otimes1_K,\\
(\operatorname{id}_A\otimes\varepsilon_K)\rho(a)
 &=\sum_i h^i(1_H)(h_i\cdot a)=1_H\cdot a=a.
\end{aligned}
\]
In particular, $\rho$ is injective. By the dual multiplication rule
in \eqref{dual-H} and the module-algebra rule for the action,
\[
\begin{aligned}
F_h\bigl(\rho(a)\rho(b)\bigr)
 &=\sum_{i,j}(h_i\cdot a)(h_j\cdot b)(h^ih^j)(h)\\
 &=\sum_{(h)}(h_{(1)}\cdot a)(h_{(2)}\cdot b)
  =h\cdot(ab)=F_h\rho(ab).
\end{aligned}
\]
Thus $\rho(ab)=\rho(a)\rho(b)$. For the involution, first note that
$S_H(S_H(h^*)^*)=h$, since $S_H^2=\operatorname{id}_H$ and $S_H$
preserves $*$. The dual $*$-operation in \eqref{dual-H} and the
action's involution rule therefore give
\[
\begin{aligned}
F_h\bigl(\rho(a)^*\bigr)
 &=\sum_i\overline{h^i(S_H(h^*))}(h_i\cdot a)^*\\
 &=\bigl(S_H(h^*)\cdot a\bigr)^*
  =S_H\bigl(S_H(h^*)^*\bigr)\cdot a^*
  =h\cdot a^*=F_h\rho(a^*).
\end{aligned}
\]
Hence $\rho$ is a unital $*$-homomorphism. Finally,
\[
\begin{aligned}
(\operatorname{id}_A\otimes\operatorname{ev}_h\otimes
 \operatorname{ev}_k)(\rho\otimes\operatorname{id}_K)\rho(a)
 &=h\cdot(k\cdot a)=(hk)\cdot a\\
 &=(\operatorname{id}_A\otimes\operatorname{ev}_h\otimes
 \operatorname{ev}_k)(\operatorname{id}_A\otimes\Delta_K)\rho(a),
\end{aligned}
\]
where the last equality uses $\Delta_K(\phi)(h\otimes k)=\phi(hk)$.
The separating slices therefore give coassociativity, so $\rho$ is
a coaction.

Conversely, let $\rho\colon A\to A\otimes K$ be a coaction and
put $h\cdot a=F_h\rho(a)$. This is bilinear, and the counit and
unitality of $\rho$ give
\[
1_H\cdot a=(\operatorname{id}_A\otimes\varepsilon_K)\rho(a)=a,
\qquad
h\cdot1_A=\operatorname{ev}_h(1_K)1_A=\varepsilon_H(h)1_A.
\]
Applying the same two-leg evaluation to coassociativity gives
$h\cdot(k\cdot a)=(hk)\cdot a$. Multiplicativity of $\rho$ and the
dual multiplication rule give
\[
h\cdot(ab)=F_h\bigl(\rho(a)\rho(b)\bigr)
 =\sum_{(h)}F_{h_{(1)}}\rho(a)\,F_{h_{(2)}}\rho(b)
 =\sum_{(h)}(h_{(1)}\cdot a)(h_{(2)}\cdot b).
\]
For $z\in A\otimes K$, the same dual involution formula used above
gives $(F_h(z))^*=F_{S_H(h^*)}(z^*)$. Since
$\rho(a)^*=\rho(a^*)$, it follows that
\[
(h\cdot a)^*=F_{S_H(h^*)}\bigl(\rho(a^*)\bigr)
             =S_H(h^*)\cdot a^*.
\]
All five action rules are therefore satisfied. Finally, the
dual-basis identities show that the two constructions are inverse:
\[
\sum_iF_{h_i}\rho(a)\otimes h^i=\rho(a),
\qquad
\sum_i h^i(h)(h_i\cdot a)=h\cdot a.
\]
\end{proof}

\subsection{Crossed products and saturation}

For a given action of $H$ on a $\Cst$-algebra $A$, the \emph{crossed product $\Cst$-algebra} is defined as follows \cite[Definition~2.6]{SP1994}:
\[
A \rtimes H : = A \otimes H
\]
with the following operations:
\[
(a \otimes h) (b \otimes k) : = \sum a(h_{(1)} \cdot b) \otimes h_{(2)} k
\]
\[
\text{and } (a \otimes h)^* := \sum h_{(1)}^* \cdot a^* \otimes h_{(2)}^*.
\]

The algebras $A$ and $H$ are embedded in $A \otimes H$ as $a \mapsto a\otimes 1_H$ and $h \mapsto 1_A \otimes h$. Using these embeddings and the above operations, for any $h\in H$ and $a\in A$, we can write $ah $ instead of $(a \otimes 1_H) (1_A \otimes h)$. Since $S^2=\operatorname{id}$, we have $\sum h_{(2)}S(h_{(1)})=\varepsilon(h)1_H$, and hence
\begin{equation*}
   ah = \sum h_{(2)} (S(h_{(1)}) \cdot a).
 \end{equation*}

We first recall the fixed-point algebra and the module notation used to
express saturation. For right modules, our inner products are linear in
the second variable.

\begin{definition}\cite[Section~4]{SP1994}\label{Facts-Hopf}
Let $H$ act on a unital $\Cst$-algebra $A$, and put $B=A\rtimes H$.
\begin{enumerate}
 \item The \emph{fixed-point algebra} is given by
 \[
  P=A^H=\{a\in A:h\cdot a=\varepsilon_H(h)a
                     \text{ for every }h\in H\}.
 \]
 \item The space $Ae_H\subseteq B$ is a left $B$-module and a right
       $P$-module, with
 \[
  (ah)(xe_H)=a(h\cdot x)e_H,\qquad (xe_H)p=xpe_H
  \quad(a,x\in A,\ h\in H,\ p\in P).
 \]
 Its $P$-valued and $B$-valued inner products are
 \begin{align*}
  \langle xe_H,ye_H\rangle_P&=e_H\cdot(x^*y),\\
  {}_B\langle xe_H,ye_H\rangle&=xe_Hy^*.
 \end{align*}
 Let $X$ be the Hilbert $P$-module completion for the first inner
 product. We identify $A$ with the dense subspace $Ae_H$ by
 $x\mapsto xe_H$.
 \item Write $\mathrm{End}(A)=\mathcal L_P(X)$ for the algebra of
       adjointable right $P$-module operators, and
 \begin{align*}
  K(A)&=\mathcal K_P(X)
   =\overline{\operatorname{span}}
       \{\theta_{\xi,\eta}:\xi,\eta\in X\},\\
  \theta_{\xi,\eta}(\zeta)&=\xi\langle\eta,\zeta\rangle_P.
 \end{align*}
 Thus $K(A)$ is the closed ideal of compact module operators in
 $\mathrm{End}(A)$.
\end{enumerate}
\end{definition}

\begin{definition}\cite[Definition~4.2]{SP1994}\label{def:saturated}
The action is \emph{saturated} if the module $X$ in
Definition~\ref{Facts-Hopf} is a $B$--$P$ imprimitivity bimodule.
Equivalently, the two inner products are full:
\[
 \overline{\operatorname{span}}\,{}_B\langle X,X\rangle=B,
 \qquad
 \overline{\operatorname{span}}\langle X,X\rangle_P=P.
\]
\end{definition}

\subsection{Conditional expectations and basic constructions}

We next recall conditional expectations, finite Watatani index,
and the basic construction. For a unital inclusion $P\subseteq A$
with $1_P=1_A$, a \emph{conditional expectation} $E:A\to P$ plays
the role of an averaging map onto the smaller algebra. More precisely,
it is a completely positive contractive map that fixes $P$ pointwise
and is $P$-bimodular:
\[
 E(p)=p,\qquad E(paq)=pE(a)q
 \quad(p,q\in P,\ a\in A).
\]
In particular, $E(1_A)=1_A$ and $E^2=E$. It is called
\emph{faithful} if $E(x^*x)=0$ implies $x=0$.
The index-finite type condition requires a finite family that reconstructs
each element of $A$ from its averaged coefficients.

\begin{definition}\cite[Definition~2.1]{OKT-2009}
Let $P\subseteq A$ be a unital inclusion of $\Cst$-algebras, and let
$E:A\to P$ be a conditional expectation.
\begin{enumerate}[label=(\alph*)]
 \item A finite family $\{(u_i,v_i)\}_{i=1}^n\subset A\times A$ is a
       \emph{quasi-basis} for $E$ if, for every $a\in A$,
       \[
        a=\sum_{i=1}^n u_iE(v_i a)
         =\sum_{i=1}^n E(a u_i)v_i.
       \]
 \item The expectation $E$ has \emph{finite Watatani index}, or is of
       \emph{index-finite type}, if it admits a quasi-basis. In this case
       its Watatani index is
       \[
        \operatorname{Index}(E)=\sum_{i=1}^n u_i v_i.
       \]
\end{enumerate}
If no quasi-basis exists, we write $\operatorname{Index}(E)=\infty$.
\end{definition}

The index is independent of the chosen quasi-basis and belongs to the
center $Z(A)$. A quasi-basis can be chosen in the form
$\{(w_i,w_i^*)\}$, so that
$\operatorname{Index}(E)=\sum_i w_iw_i^*\geq0$. When $A$ is simple,
the index is therefore a positive scalar multiple of $1_A$.
Being of index-finite type also ensures that $E$ is faithful; see \cite[Remark~2.2]{OKT-2009}.
The index is moreover invertible, as recalled in
\cite[Remark~2.3(4)]{OKT-2009}.

We next describe the basic construction, following
\cite[Section~2.1.2]{OKT-2009}. Assume that $E$ has finite Watatani
index. Regard $A$ as a right $P$-module, with $P$-valued inner product
\[
 \langle x,y\rangle_P=E(x^*y),\qquad x,y\in A.
\]
Let $\mathcal E$ be its Hilbert $P$-module completion and let
$\eta:A\to\mathcal E$ be the natural inclusion. Faithfulness of $E$
ensures that $\eta$ is injective. Write $\mathcal L_P(\mathcal E)$ for
the $\Cst$-algebra of adjointable right $P$-linear operators on
$\mathcal E$. Here an adjoint $T^*$ satisfies
\[
 \langle T\xi,\zeta\rangle_P
   =\langle\xi,T^*\zeta\rangle_P
 \quad(\xi,\zeta\in\mathcal E).
\]
Left multiplication by $A$ and averaging by $E$ give operators
\[
 \lambda(a)\eta(x)=\eta(ax),\qquad
 e_P\eta(x)=\eta(E(x)).
\]
The map $\lambda:A\to\mathcal L_P(\mathcal E)$ is an injective
$*$-homomorphism, and $e_P$ is a projection, called the
\emph{Jones projection}. It records the expectation through the
relation
\[
 e_P\lambda(a)e_P=\lambda(E(a))e_P,\qquad a\in A.
\]
The \emph{reduced $\Cst$-basic construction} is
\[
 C_r^*\langle A,e_P\rangle
  =\overline{\operatorname{span}}
       \{\lambda(a)e_P\lambda(b):a,b\in A\}
  \subseteq\mathcal L_P(\mathcal E).
\]
We identify $A$ with $\lambda(A)$ and write $ae_Pb$ for
$\lambda(a)e_P\lambda(b)$. Thus the basic construction is an algebra
encoding both $A$ and its averaging projection. The associated
\emph{dual conditional expectation} maps this algebra back onto $A$:
\[
 \widehat E:C_r^*\langle A,e_P\rangle\to A,\qquad
 \widehat E(ae_Pb)=\operatorname{Index}(E)^{-1}ab;
\]
see \cite[Remark~2.3(4)]{OKT-2009}.

We now return to the Hopf algebra action. We put $P=A^H$, $B=A\rtimes H$,
and $N=\dim H$. Let $\tau_H$ be the normalized Haar trace on $H$;
then $\tau_H(e_H)=N^{-1}$ by \cite[Theorem~2.2]{SP1994}.
The two canonical averaging maps of Szyma\'nski and Peligrad are
\[
 \begin{aligned}
  E:A\longrightarrow P,&\qquad E(a)=e_H\cdot a,\\
  F:B\longrightarrow A,&\qquad F(ah)=\tau_H(h)a
                      \quad(a\in A,\ h\in H),
 \end{aligned}
\]
where the formula for $F$ is extended linearly from elementary
products $ah$. By \cite[Proposition~2.12]{SP1994}, $E$ is a faithful
conditional expectation and satisfies
\[
 E((h\cdot a)b)=E(a(S_H(h)\cdot b))
 \quad(h\in H,\ a,b\in A).
\]
By \cite[Proposition~2.8]{SP1994}, $F$ is also a faithful conditional
expectation. These assertions do not require saturation. The following
proposition explains how the two expectations are related when the
action is saturated and note that no simplicity assumption on $A$ or $B$ is needed.

\begin{proposition}\label{prop:crossed-product-basic-construction}
Let $H$ act on a unital $\Cst$-algebra $A$, and suppose that the
action is saturated. Let $E:A\to P$ and $F:B\to A$ be the canonical expectations above.
Then $\operatorname{Index}(E)=N1_A$, and there is a unital
$*$-isomorphism fixing $A$,
\[
 \Phi:C_r^*\langle A,e_P\rangle\longrightarrow B,
 \qquad \Phi(ae_Pb)=ae_Hb.
\]
Under this identification, $F$ is the dual expectation of $E$:
\[
 F\circ\Phi=\widehat E.
\]
\end{proposition}

\begin{proof}
By \cite[Proposition~4.5]{SP1994}, $E$ has index $N1_A$, and
there is a $*$-isomorphism
\[
 B\longrightarrow K(A)=\mathrm{End}(A),\qquad
 ae_Hb\longmapsto \lambda(a)e_P\lambda(b).
\]
Here the module in Definition~\ref{Facts-Hopf} is identified with the
Hilbert $P$-module of the basic construction by $a\mapsto ae_H$:
its inner product is $e_H\cdot(a^*b)=E(a^*b)$. The range of the displayed
map is therefore precisely $C_r^*\langle A,e_P\rangle$. Its inverse
is the asserted map $\Phi$. The map is unital and fixes $A$, since
the action of $a\in A$ on the module is left multiplication
$\lambda(a)$.

Write $\Delta_H(e_H)=\sum_{(e_H)}e_{H(1)}\otimes e_{H(2)}$.
Crossed-product multiplication and the identity
$he_H=\varepsilon_H(h)e_H$ give the Jones relation directly:
\begin{align*}
 e_Hae_H
 &=\sum_{(e_H)}(e_{H(1)}\cdot a)e_{H(2)}e_H\\
 &=\sum_{(e_H)}(e_{H(1)}\cdot a)
                         \varepsilon_H(e_{H(2)})e_H\\
 &=(e_H\cdot a)e_H=E(a)e_H.
\end{align*}
The third equality uses the counit identity
$\sum_{(e_H)}e_{H(1)}\varepsilon_H(e_{H(2)})=e_H$
and linearity of the action. Thus $e_H$ satisfies the same Jones
relation as $e_P$, and $\Phi(e_P)=e_H$.

Finally, $F$ is $A$-bimodular and
$F(e_H)=\tau_H(e_H)1_A=N^{-1}1_A$. Hence, for every $a,b\in A$,
\begin{align*}
 F\bigl(\Phi(ae_Pb)\bigr)
 &=F(ae_Hb)\\
 &=aF(e_H)b\\
 &=N^{-1}ab\\
 &=\widehat E(ae_Pb).
\end{align*}
The last equality is the defining formula for the dual expectation,
since $\operatorname{Index}(E)=N1_A$. The linear span of the elements
$ae_Pb$ is dense in $C_r^*\langle A,e_P\rangle$, and both maps are
bounded. Their equality therefore extends to the entire basic
construction.
\end{proof}

\section{Rokhlin property, approximate representability, and duality}\label{sec:main}

Let $H\curvearrowright A$ be an action on a unital $\Cst$-algebra.
We retain the notation of Section~\ref{sec:prelim} and will use throughout;
$P=A^H$ and $E:A\to P$ denote the fixed-point algebra and the
canonical conditional expectation, respectively.

The \emph{sequence algebra} of $A$ is
\[
 c_0(A)=\{(a_n)\in\ell^\infty(\bN,A):\|a_n\|\to0\},
 \qquad
 A^\infty=\ell^\infty(\bN,A)/c_0(A).
\]
We identify $A$ with its image in $A^\infty$ under constant sequences.
The action of $H$ induces the coordinatewise action
\[
 h^\infty\cdot[(a_n)]=[(h\cdot a_n)]
 \qquad(h\in H).
\]
For the associated coaction $\rho:A\to A\otimes K$, write
$\rho^\infty:A^\infty\to A^\infty\otimes K$
for the induced coaction, also defined coordinatewise.

\subsection{The Rokhlin property}

In \cite{OKT-2009}, the Rokhlin property for actions of $H$ was originally defined by the existence of an equivariant representation with values in the central sequence algebra
\[
\pi:K \longrightarrow A^\infty\cap A'.
\]
However, \cite{KT-2015} shows that this strong requirement forces the dual of \(H\) to be commutative. Therefore, we propose the following modification.
\begin{definition}\label{def:RP-Hopf}
  An action of a finite-dimensional Hopf $\Cst$-algebra $H$ on $A$ is said to have the Rokhlin property if
there exists a unital $*$-homomorphism, called a \emph{Rokhlin homomorphism},
$\pi: K \to A^{\infty}$ such that
\begin{enumerate}
\item \[h^{\infty} \cdot \pi(\varphi)= \sum_{(\varphi)}\varphi_{(2)} (h) \pi(\varphi_{(1)}) \quad\text{for } h \in H \text{ and } \varphi \in K;\]
\item \[ \pi(e_K) \in A^{\infty}\cap A', \] 
\end{enumerate} 
where $e_K$ is the distinguished projection in $K$. 
\end{definition}

\begin{lemma}\label{lem:reconstruct-rokhlin-pi}
Let $H$ be a finite-dimensional Hopf $\Cst$-algebra with $N=\dim H$. Let $e=e_H$ and $e_K$ be the distinguished
projections of $H$ and $K$, respectively. Suppose that
$\pi:K\to A^\infty$ satisfies the Rokhlin conditions
in Definition~\ref{def:RP-Hopf}. If $p=\pi(e_K)$, then
\[
 \pi(\phi)
 =
 \left(N\sum_{(e)}\phi(e_{(1)})S(e_{(2)})\right)\cdot p
 \qquad (\phi\in K).
\]
Here $S$ is the antipode of $H$, and $\cdot$ denotes the induced action
on $A^\infty$.
\end{lemma}

\begin{proof}
Write $\tau=\tau_H=e_K$, where the equality is as functionals on $H$:
Haar invariance gives
$\phi\tau=\tau\phi=\phi(1_H)\tau=\varepsilon^\circ(\phi)\tau$ for
$\phi\in K$, and $\varepsilon^\circ(\tau)=\tau(1_H)=1$, which determine
$e_K$ uniquely.
Set
\[ q_\phi=N\sum_{(e)}\phi(e_{(1)})S(e_{(2)}). \]

For $x,y\in H$, apply
\[
 (\tau\otimes\operatorname{id}_H)\Delta_H(z)=\tau(z)1_H
\]
to $z=xy_{(1)}$, keeping the remaining tensor leg $y_{(2)}$.
Since $\Delta_H$ is multiplicative, this gives
\[
 \sum_{(y)}\tau(xy_{(1)})1_H\otimes y_{(2)}
 =
 \sum_{(x),(y)}
 \tau(x_{(1)}y_{(1)})x_{(2)}y_{(2)}\otimes y_{(3)}.
\]
Applying $m_H(\operatorname{id}_H\otimes S)$ yields
\begin{align*}
 \sum_{(y)}\tau(xy_{(1)})S(y_{(2)})
 &=
 \sum_{(x),(y)}
 \tau(x_{(1)}y_{(1)})x_{(2)}y_{(2)}S(y_{(3)})\\
 &=
 \sum_{(x),(y)}
 \tau(x_{(1)}y_{(1)})x_{(2)}\varepsilon_H(y_{(2)})\\
 &=
 \sum_{(x)}\tau(x_{(1)}y)x_{(2)}.
\end{align*}
For the second equality, coassociativity lets us read the three
$y$-legs as
\[
 (\operatorname{id}_H\otimes\Delta_H)\Delta_H(y).
\]
The antipode identity on the last two legs then gives
\[
 \sum_{(y)}y_{(1)}\otimes y_{(2)}S(y_{(3)})
 =
 \sum_{(y)}y_{(1)}\otimes\varepsilon_H(y_{(2)})1_H.
\]
The last equality uses the counit identity
$\sum_{(y)}y_{(1)}\varepsilon_H(y_{(2)})=y$. Taking $y=e$, we obtain
\begin{align*} 
 \sum_{(e)}\tau(xe_{(1)})S(e_{(2)})
 &=
 \sum_{(x)}\tau(x_{(1)}e)x_{(2)}\notag\\
 &=
 \frac{1}{N}\sum_{(x)}\varepsilon_H(x_{(1)})x_{(2)}
 =\frac{x}{N}.
 \end{align*}
Moreover, for all $y \in H$, we have
\[ 
\begin{split}
 \tau \left(y N \sum_{(e)} e_{(1)}\tau(xS(e_{(2)}))\right)&=N \sum_{(e)} \tau (xS(e_{(2)}))\tau(ye_{(1)})\\
 &= N \tau \left(x \sum_{(e)} \tau(ye_{(1)})S(e_{(2)})\right)\\
 &= \tau(xy)=\tau(yx), 
\end{split}
\] 
and hence
\[ N \sum_{(e)} e_{(1)}\tau(xS(e_{(2)})) =x. \]

Finally, since $p=\pi(\tau)$, condition~(1) in Definition~\ref{def:RP-Hopf} gives
\[
 q_\phi\cdot p
 =
 \pi\!\left(\sum_{(\tau)}
 \tau_{(2)}(q_\phi)\tau_{(1)}\right)
 .
\]

We claim that 

\[
 \sum_{(\tau)}\tau_{(1)}(x)\tau_{(2)}(q_\phi)
 =\phi(x).
\]
Indeed, 
\[ 
\begin{split} 
\sum_{(\tau)}\tau_{(1)}(x)\tau_{(2)}(q_\phi)&=\tau(q_{\phi}x)\\
&=N \sum_{(e)} \phi(e_{(1)})\tau(S(e_{(2)})x) \\
&= \phi \left(N \sum_{(e)} e_{(1)} \tau(S(e_{(2)})x)\right)\\
&= \phi(x)
\end{split}
\]
Hence that functional is $\phi$, and $q_\phi\cdot p=\pi(\phi)$.
Substituting the definition of $q_\phi$ proves the formula.
\end{proof}

\begin{definition}\cite[Definition 3.1]{OKT-2009}\label{def:RP-incl}
 Let $E:A\to P$ be a conditional expectation of index-finite type between unital $\Cst$-algebras. Then $E$ is said to have the \emph{Rokhlin property} if there exists a projection $e \in A^{\infty}\cap A'$ such that 
\[
E^{\infty}(e) = (\operatorname{Index}E)^{-1} \cdot 1
\]
and the map $A \ni x  \mapsto xe$ is injective.
\end{definition}

For the associated coaction $\rho\colon A\to A\otimes K$,
Kodaka--Teruya's \cite[Theorem~5.12 and Corollary~6.4]{KT-2015}
already characterize the Rokhlin property by a projection
$p\in A^\infty\cap A'$ with $E^\infty(p)=(\dim H)^{-1}1_A$;
the latter result does not assume saturation. The following theorem
relates that established criterion and our direct representation
formulation to the Rokhlin property of the
fixed-point expectation and provides the Rokhlin part of the
inclusion-theoretic comparison, complementary to
Theorem~\ref{thm:main} for approximate representability.

\begin{theorem}\label{thm:Hopf-incl}
Let $H$ be a finite-dimensional Hopf $\Cst$-algebra with $N=\dim H$, and let
$H\curvearrowright A$ be a saturated action on a unital simple
$\Cst$-algebra. Let $E\colon A\to A^H$ be the canonical conditional
expectation from Proposition~\ref{prop:crossed-product-basic-construction}. Then
$H\curvearrowright A$ has the Rokhlin property if and only if $E$
has the Rokhlin property.
\end{theorem}
\begin{proof}
Write $j:H\to B$ for the canonical embedding and put $f=j(e_H)$.
In the sense of Szyma\'nski--Peligrad, saturation means that the
$B$-valued inner products of $Af$ span a dense subspace of $B$;
these inner products are ${}_B\langle af,bf\rangle=afb^*$.
Equivalently,
\[
 B=\overline{\operatorname{span}}\{afb:a,b\in A\}.
\]
Under this hypothesis, $B$ identifies with the basic construction for
$P\subset A$, with $f$ corresponding to the Jones projection, and
$\operatorname{Index}(E)=N1_A$; see
\cite[Definition~4.2 and Proposition~4.5]{SP1994}.

\smallskip
\noindent\emph{From the action to the expectation.}
Let $\pi:K\to A^\infty$ satisfy the Rokhlin conditions in
Definition~\ref{def:RP-Hopf}, and put $p=\pi(e_K)$.
Then $p$ is a projection in $A^\infty\cap A'$.
For the associated coaction $\rho:A\to A\otimes K$ we have
$\rho(a)=\sum_i(h_i\cdot a)\otimes h^i$. The Haar trace on $K$ is
$\tau_K(\phi)=\phi(e_H)$, so the dual-basis identity gives
\begin{align*}
 (\operatorname{id}_A\otimes\tau_K)\rho(a)
 &=\sum_i \tau_K(h^i)(h_i\cdot a)\\
 &=\left(\sum_i h^i(e_H)h_i\right)\cdot a\\
 &=e_H\cdot a=E(a).
\end{align*}
Applying the same identity coordinatewise gives
$E^\infty=(\operatorname{id}_{A^\infty}\otimes\tau_K)\circ\rho^\infty$.
Equivariance gives
$\rho^\infty(\pi(\phi))=\sum_{(\phi)}\pi(\phi_{(1)})\otimes\phi_{(2)}$.
Thus
\begin{align*}
 E^\infty(\pi(\phi))
 &=\sum_{(\phi)}\pi(\phi_{(1)})\tau_K(\phi_{(2)})\\
 &=\pi\!\left(\sum_{(\phi)}\phi_{(1)}\tau_K(\phi_{(2)})\right)\\
 &=\pi(\tau_K(\phi)1_K)=\tau_K(\phi)1_A,
\end{align*}
where the last line uses Haar invariance and unitality of $\pi$.
Since $e_K=\tau_H$ as a functional on $H$, we obtain
\[
 E^\infty(p)=\tau_K(e_K)1_A
 =e_K(e_H)1_A=\tau_H(e_H)1_A=N^{-1}1_A.
\]
In particular, $p\ne0$. Centrality of $p$ relative to $A$ gives
$(ap)(bp)=abp$ and $(ap)^*=a^*p$ for $a,b\in A$. Therefore
$a\mapsto ap$ is a nonzero $*$-homomorphism, and its kernel is zero
by simplicity of $A$. Hence $p$ is a Rokhlin projection for $E$.

\smallskip
\noindent\emph{From the expectation to the action.}
Suppose $E$ has a Rokhlin projection $p\in A^\infty\cap A'$.
Then $E^\infty(p)=N^{-1}1_A$. By
Proposition~\ref{prop:crossed-product-basic-construction}, the canonical
expectation $F:B\to A$ is the dual expectation of $E$; in particular,
$F(aj(h))=a\tau_H(h)$ and $F(f)=N^{-1}1_A$.
We work in $B^\infty$, identifying $A^\infty$ with its canonical
image. The finite-index inequality for $F$ passes to sequences,
so $F^\infty$ is faithful. From $fxf=E^\infty(x)f$ for $x\in A^\infty$
we get $fpf=f/N$. Thus $q=Npfp$ is a projection below $p$, since
$q^*=q$ and $q^2=N^2p(fpf)p=q$. Bimodularity gives
$F^\infty(q)=p$, so faithfulness applied to $p-q$ gives $q=p$.
Consequently $pfp=p/N$. For $a,b\in A$,
\[
 p(afb)p=a(pfp)b=N^{-1}abp=F(afb)p.
\]
Saturation and continuity extend this identity to every element of
$B$. In particular, $pj(h)p=\tau_H(h)p$ for $h\in H$.

Write $\tau=\tau_H$ and $S=S_H$. The crossed-product covariance
formula, applied coordinatewise, now gives
\begin{align*}
 p(h\cdot p)
 &=\sum_{(h)}pj(h_{(1)})p\,j(S(h_{(2)}))\\
 &=\sum_{(h)}\tau(h_{(1)})p\,j(S(h_{(2)}))
 =\tau(h)p.
\end{align*}
For the last equality apply $S$ to the Haar invariance identity
$\sum_{(h)}\tau(h_{(1)})h_{(2)}=\tau(h)1_H$.
The induced action on $A^\infty$ satisfies the module-algebra rule
\[
 h\cdot(xy)=\sum_{(h)}(h_{(1)}\cdot x)(h_{(2)}\cdot y)
 \qquad (x,y\in A^\infty).
\]
Fix $h,k\in H$. We use Sweedler notation
\[
 \Delta_H(h)=\sum_{(h)}h_{(1)}\otimes h_{(2)},
 \qquad
 (\operatorname{id}_H\otimes\Delta_H)\Delta_H(h)
 =\sum_{(h)}h_{(1)}\otimes h_{(2)}\otimes h_{(3)}.
\]
By coassociativity, the same three-leg sum also represents
$(\Delta_H\otimes\operatorname{id}_H)\Delta_H(h)$.
Set
\[
 T=\sum_{(h)}h_{(1)}\cdot
       \bigl(p((S(h_{(2)})k)\cdot p)\bigr).
\]
To expand $T$, apply the module-algebra rule to the outer action
of $h_{(1)}$ on the product of $p$ and $(S(h_{(2)})k)\cdot p$.
This splits the first coproduct leg, so the resulting three legs
come from $(\Delta_H\otimes\operatorname{id}_H)\Delta_H(h)$.
Then use $a\cdot(b\cdot x)=(ab)\cdot x$. Thus
\begin{align*}
 T&=\sum_{(h)}(h_{(1)}\cdot p)
       \bigl(h_{(2)}\cdot((S(h_{(3)})k)\cdot p)\bigr)\\
  &=\sum_{(h)}(h_{(1)}\cdot p)
       \bigl((h_{(2)}S(h_{(3)})k)\cdot p\bigr).
\end{align*}
To apply the antipode axiom, use coassociativity to read this same
sum as
\[
 (\operatorname{id}_H\otimes\Delta_H)\Delta_H(h).
\]
The first leg is unchanged, and the last two legs now come from splitting the
second leg of $\Delta_H(h)$. The ordinary antipode axiom
\(
 m_H(\operatorname{id}_H\otimes S)\Delta_H(a)
 =\varepsilon_H(a)1_H
\)
applies to that inner sum, with $a$ the unsplit second leg. By
linearity of the action and the counit identity
$\sum_{(h)}h_{(1)}\varepsilon_H(h_{(2)})=h$, we obtain
\begin{align*}
 T&=\sum_{(h)}\varepsilon_H(h_{(2)})(h_{(1)}\cdot p)(k\cdot p)\\
  &=\left(\left[\sum_{(h)}h_{(1)}\varepsilon_H(h_{(2)})\right]\cdot p\right)
       (k\cdot p)\\
  &=(h\cdot p)(k\cdot p).
\end{align*}
On the other hand, use $p(\ell\cdot p)=\tau(\ell)p$ with
$\ell=S(h_{(2)})k$ directly in the defining sum for $T$. Comparing
this evaluation with the preceding one gives
\[
 (h\cdot p)(k\cdot p)
 =\sum_{(h)}\tau(S(h_{(2)})k)(h_{(1)}\cdot p).
\]

For $\phi\in K$, set
$q_\phi=N\sum_{(e_H)}\phi(e_{H(1)})S(e_{H(2)})$
and define $\pi_p(\phi)=q_\phi\cdot p$. The proof of Lemma~\ref{lem:reconstruct-rokhlin-pi}
gives $\tau(xq_\phi)=\phi(x)$ for $x\in H$, independently
of the existence of a Rokhlin homomorphism. Faithfulness of $\tau$
makes this pairing nondegenerate, so it determines $q_\phi$ uniquely.
For $r=q_\phi$ and $t=q_\chi$, the product rule just proved gives
$\pi_p(\phi)\pi_p(\chi)=v\cdot p$, where
$v=\sum_{(r)}\tau(S(r_{(2)})t)r_{(1)}$. Using the Haar identity
$\sum_{(r)}\tau(xr_{(1)})S(r_{(2)})
 =\sum_{(x)}\tau(x_{(1)}r)x_{(2)}$ from the preceding proof, we have
\begin{align*}
 \tau(xv)
 &=\sum_{(r)}\tau(xr_{(1)})\tau(S(r_{(2)})t)\\
 &=\sum_{(x)}\tau(x_{(1)}r)\tau(x_{(2)}t)\\
 &=\sum_{(x)}\phi(x_{(1)})\chi(x_{(2)})=(\phi\chi)(x).
\end{align*}
Thus $v=q_{\phi\chi}$ and $\pi_p$ is multiplicative.
For the involution, $\tau\circ S=\tau$ and $S^2=\operatorname{id}_H$
give
\[
 \phi^*(x)=\overline{\phi(S(x)^*)}
 =\overline{\tau(S(x)^*q_\phi)}
 =\tau(q_\phi^*S(x))=\tau(xS(q_\phi^*)).
\]
Hence $q_{\phi^*}=S(q_\phi^*)$ and the $*$-rule for the action gives
$\pi_p(\phi)^*=\pi_p(\phi^*)$. The defining pairing gives
$q_{1_K}=Ne_H$ and $q_{e_K}=1_H$, since
$\tau(xNe_H)=\varepsilon_H(x)$ and $e_K(x)=\tau(x)$. Therefore
\[
 \pi_p(1_K)=N E^\infty(p)=1_A,\qquad \pi_p(e_K)=p.
\]
Finally, if $(R_h\phi)(x)=\phi(xh)$, then
$\tau(x hq_\phi)=\phi(xh)$ gives $q_{R_h\phi}=hq_\phi$. Consequently
\[
 h\cdot\pi_p(\phi)=\pi_p(R_h\phi)
 =\sum_{(\phi)}\phi_{(2)}(h)\pi_p(\phi_{(1)}),
\]
where the last equality follows from
$\Delta_K(\phi)(x\otimes h)=\phi(xh)$. Thus $\pi_p$ is equivariant,
and $\pi_p(e_K)=p\in A^\infty\cap A'$ proves the Rokhlin property of
the action.
\end{proof}

For a finite-dimensional Hopf $\Cst$-algebra $H$, write
$H^{\mathrm{cop}}=(H,\Delta_H^{\mathrm{cop}})$ for its co-opposite
Hopf $\Cst$-algebra. It has the same underlying algebra, involution,
and counit as $H$, with coproduct
\[
 \Delta_H^{\mathrm{cop}}=\Sigma_{H,H}\circ\Delta_H,
\]
where $\Sigma_{X,Y}:X\otimes Y\to Y\otimes X$ is the tensor flip,
$\Sigma_{X,Y}(x\otimes y)=y\otimes x$. The antipode of
$H^{\mathrm{cop}}$ is $S_H^{-1}$.

For the spatial Rokhlin comparison, let
$\mathbb G$ be the finite quantum group with
\[
 C(\mathbb G)=K,\qquad
 \Delta_{\mathbb G}=\Delta_K^{\mathrm{cop}}
                  =\Sigma_{K,K}\circ\Delta_K.
\]
Thus the Hopf $\Cst$-algebra of functions on $\mathbb G$ is
$K^{\mathrm{cop}}$. If $\rho:A\to A\otimes K$ is the coaction
associated with $H\curvearrowright A$, define
\[
 \lambda:=\Sigma_{A,K}\circ\rho:A\longrightarrow C(\mathbb G)\otimes A,
 \qquad
 \lambda(a)=\sum_i h^i\otimes(h_i\cdot a).
\]
It is easily verified, by flipping the tensor factors in the
coassociativity and counit identities for $\rho$, that
\begin{align*}
 (\Delta_{\mathbb G}\otimes\operatorname{id}_A)\lambda
  &=(\operatorname{id}_K\otimes\lambda)\lambda,\\
 (\varepsilon_K\otimes\operatorname{id}_A)\lambda
  &=\operatorname{id}_A.
\end{align*}
Thus $\lambda$ gives the action of $\mathbb G$ on $A$ in the
convention of Barlak--Szab\'o--Voigt.

\Needspace{8\baselineskip}
\begin{theorem}\label{thm:RP-SRP}
With the notation above, let $A$ be separable and unital.
If the action $H\curvearrowright A$ is saturated, then it has the
Rokhlin property if and only if the action $\lambda$ of $\mathbb G$
on $A$ has the spatial Rokhlin property in the sense of
Barlak--Szab\'o--Voigt. The implication from the spatial Rokhlin
property to the Rokhlin property of $H\curvearrowright A$ holds
without assuming saturation.
\end{theorem}

\begin{proof}
Write $K=C(\mathbb G)=H^\circ$, $N=\dim H$ and
$\delta=\Delta_{\mathbb G}=\Delta_K^{\mathrm{cop}}$. The symbol $K$
by itself carries the original Hopf structure; we use $\delta$ for
the coproduct of $\mathbb G$. Its antipode is
$s=S_K^{-1}=S_K$, since $S_K^2=\operatorname{id}_K$.
The preliminary identifications below hold for every action;
saturation will be used only in the forward implication.

By \cite[Definitions~3.1 and~4.1]{BSV2017}, the spatial Rokhlin property
means that the canonical embedding
\[
 j_A:A\longrightarrow D:=C(\mathbb G)\boxtimes A,
 \qquad j_A(a)=1\boxtimes a,
\]
with $D$ carrying the coaction $\Omega=\delta\boxtimes\lambda$,
admits an equivariant sequential splitting. That is, there is a
$*$-homomorphism $\psi:D\to A^\infty$ such that
\[
 \psi(j_A(a))=a,\qquad
 \lambda^\infty\psi=(\operatorname{id}_K\otimes\psi)\Omega,
\]
where $A$ is identified with its constant sequences. Since $K$ is a
finite direct sum of matrix algebras, the canonical map
$K\otimes A^\infty\to(K\otimes A)^\infty$ is an isomorphism
(entrywise on each matrix summand). Thus $\lambda^\infty$ is given
on all of $A^\infty$ by applying $\lambda$ coordinatewise.

The finite quantum group $\mathbb G$ is regular and coexact (the
latter follows from exactness of the finite-dimensional algebra $K$).
For completeness, $\lambda$ is unital and injective, and it satisfies
the coassociativity and counit identities above. To see the density
condition, write $\lambda(a)=\sum_i h^i\otimes a_i$ with
$a_i=h_i\cdot a$, and set
\[
 \lambda(a_i)=\sum_j h^j\otimes a_{ij},\qquad
 \delta(h^r)=\sum_{i,j}d^r_{ij}h^i\otimes h^j.
\]
Comparing the $h^i\otimes h^j$ coefficients of
$(\operatorname{id}_K\otimes\lambda)\lambda(a)$ and
$(\delta\otimes\operatorname{id}_A)\lambda(a)$ gives
$a_{ij}=\sum_r d^r_{ij}a_r$. Hence, for every $\phi\in K$,
\begin{align*}
 \sum_i(\phi s(h^i)\otimes1)\lambda(a_i)
 &=\sum_{i,j}\phi s(h^i)h^j\otimes a_{ij}\\
 &=\sum_r\phi m_K(s\otimes\operatorname{id}_K)\delta(h^r)\otimes a_r\\
 &=\sum_r\varepsilon_K(h^r)\phi\otimes a_r=\phi\otimes a.
\end{align*}
Here the third equality is the antipode identity for $\delta$ and
the last one is the counit identity for $\lambda$. Taking adjoints
gives the density condition with the factors in the opposite order
as well.

By \cite[Lemma~1.24]{BSV2017}, there is an equivariant $*$-isomorphism
\[
 T_\lambda:(D,\Omega)\longrightarrow (K\otimes A,\Gamma),
 \qquad \Gamma=\delta\otimes\operatorname{id}_A,
\]
with
\[
 T_\lambda(\phi\boxtimes1)=\phi\otimes1,
 \qquad T_\lambda(1\boxtimes a)=\lambda(a).
\]
In particular, in ordinary tensor coordinates the embedding to be
split is $\lambda$, not $a\mapsto1\otimes a$. Equivariance of
$T_\lambda$ gives $\Gamma T_\lambda=(\operatorname{id}_K\otimes
T_\lambda)\Omega$ and thus
$\Omega T_\lambda^{-1}=(\operatorname{id}_K\otimes T_\lambda^{-1})\Gamma$.

\smallskip
\noindent\emph{From Hopf Rokhlin property to spatial Rokhlin property.}
Suppose $H\curvearrowright A$ is saturated and has the Rokhlin property, with a
unital $*$-homomorphism $\pi:K\to A^\infty$ as in
Definition~\ref{def:RP-Hopf}. Put $p=\pi(e_K)\in A^\infty\cap A'$.
Flipping the coaction form of equivariance of $\pi$ gives
\[
 \lambda^\infty(\pi(\phi))
 =\sum_{(\phi)}\phi_{(2)}\otimes\pi(\phi_{(1)})
 =(\operatorname{id}_K\otimes\pi)\delta(\phi),
\]
where the Sweedler indices on $\phi$ are for $\Delta_K$.
The averaging calculation in the forward direction of
Theorem~\ref{thm:Hopf-incl} uses only equivariance and unitality
of $\pi$ and gives
\[
 E^\infty(p)=\tau_K(e_K)1_A=\tau_H(e_H)1_A=N^{-1}1_A.
\]

Following the construction in \cite[Lemma~5.1 and the proof of
Theorem~5.7]{OKT-2009}, define
\[
 \beta(a)=N E^\infty(ap)\in P^\infty,\qquad a\in A.
\]
We verify directly that $\beta$ is a unital $*$-homomorphism into
$\{p\}'\cap P^\infty$; we do not need injectivity of $\beta$ or
assume that $E$ has the full Rokhlin property.

By saturation the basic construction is $B=A\rtimes H$. Let
$f=j(e_H)$ be its Jones projection and $F:B\to A$ its dual
conditional expectation. As in the proof of
Theorem~\ref{thm:Hopf-incl}, $F(f)=N^{-1}1_A$ and $F^\infty$ is
faithful by the finite-index order bound for $F$. The Jones
relation on sequences gives
$fpf=E^\infty(p)f=N^{-1}f$. Hence $q=Npfp$ satisfies $q^*=q$,
$q^2=N^2p(fpf)p=q$ and $pq=qp=q$. Thus $p-q\geq0$, and
bimodularity gives $F^\infty(q)=NpF(f)p=p$. Faithfulness applied
to $p-q$ gives $q=p$, or $pfp=N^{-1}p$.

The identity $\beta(a)p=ap$ follows exactly from this compression:
using $p=Npfp$, the Jones relation, and bimodularity, in that order,
\begin{align*}
 ap&=N F^\infty(fap)\\
   &=N^2 F^\infty(fapfp)\\
   &=N^2 F^\infty(E^\infty(ap)fp)\\
   &=N E^\infty(ap)p=\beta(a)p.
\end{align*}
For $y\in P^\infty$, $E^\infty(yp)=yE^\infty(p)=N^{-1}y$.
Thus $yp=0$ implies $y=0$. In particular, $\beta(a)$ is the unique
$y\in P^\infty$ such that $yp=ap$.
The formula for $\beta$ gives $\beta(1)=1$, and since $p$ commutes
with $A$, it gives
$\beta(a)^*=N E^\infty(pa^*)=N E^\infty(a^*p)=\beta(a^*)$.
Applying $\beta(a^*)p=a^*p$ and taking adjoints yields
$p\beta(a)=pa=ap=\beta(a)p$. Finally,
\[
 \beta(a)\beta(b)p=\beta(a)bp=\beta(a)pb=apb=abp=\beta(ab)p.
\]
The difference $\beta(a)\beta(b)-\beta(ab)$ lies in $P^\infty$,
so the injectivity of $y\mapsto yp$ on $P^\infty$ proves
multiplicativity. This verifies the required map
$\beta:A\to\{p\}'\cap P^\infty$.

We next prove that $\beta(A)$ commutes with $\pi(K)$. Since each
entry of a representative in $P^\infty$ is fixed, we have
$h\cdot\beta(a)=\varepsilon_H(h)\beta(a)$; equivalently,
$\lambda^\infty(\beta(a))=1_K\otimes\beta(a)$. Apply $h\in H$ to
$p\beta(a)=\beta(a)p$. Using the product rule for the action followed
by the counit identities on $H$, the two sides become
\begin{align*}
 h\cdot(p\beta(a))
 &=\sum_{(h)}(h_{(1)}\cdot p)\varepsilon_H(h_{(2)})\beta(a)
   =(h\cdot p)\beta(a),\\
 h\cdot(\beta(a)p)
 &=\sum_{(h)}\varepsilon_H(h_{(1)})\beta(a)(h_{(2)}\cdot p)
   =\beta(a)(h\cdot p).
\end{align*}
Now Lemma~\ref{lem:reconstruct-rokhlin-pi} gives
\[
 \pi(\phi)=q_\phi\cdot p,
 \qquad q_\phi=N\sum_{(e_H)}\phi(e_{H(1)})S_H(e_{H(2)}).
\]
Taking $h=q_\phi$ shows that $\pi(\phi)\beta(a)=\beta(a)\pi(\phi)$.

We also need the following expression for $\beta$:
\begin{equation}\label{eq:spatial-rokhlin-beta}
 \beta(a)=\sum_i\pi(s(h^i))a_i.
\end{equation}
To prove this, use $s(h^i)(x)=h^i(S_H(x))$ and the dual-basis
identity $\sum_i h^i(x)a_i=x\cdot a$. The antipode reverses
the coproduct:
\[
 \Delta_H(S_H(e_H))
 =\sum_{(e_H)}S_H(e_{H(2)})\otimes S_H(e_{H(1)}).
\]
Thus the product rule for the action combines the two factors
with these coproduct legs into $S_H(e_H)\cdot(pa)$. Since
$S_H(e_H)=e_H$ and $pa=ap$, the formula for $\pi$ gives
\begin{align*}
 \sum_i\pi(s(h^i))a_i
 &=N\sum_{(e_H),i}h^i(S_H(e_{H(1)}))
     (S_H(e_{H(2)})\cdot p)a_i\\
 &=N\sum_{(e_H)}(S_H(e_{H(2)})\cdot p)
                         (S_H(e_{H(1)})\cdot a)\\
 &=N S_H(e_H)\cdot(pa)\\
 &=N E^\infty(ap)=\beta(a).
\end{align*}

Define the ordinary tensor-product map
$L:K\otimes A\to A^\infty$ by $L(\phi\otimes a)=\pi(\phi)\beta(a)$.
The two ranges commute, so on elementary tensors we have
\begin{align*}
 L(\phi\otimes a)L(\chi\otimes b)
 &=\pi(\phi)\pi(\chi)\beta(a)\beta(b)
   =L(\phi\chi\otimes ab),\\
 L(\phi\otimes a)^*
 &=\beta(a^*)\pi(\phi^*)=L(\phi^*\otimes a^*).
\end{align*}
Both maps are unital. Their commuting ranges give a unital
$*$-homomorphism on the maximal tensor product, which equals
the minimal tensor product because $K$ is finite-dimensional.
In particular, $L$ is well-defined on $K\otimes A$.

Using \eqref{eq:spatial-rokhlin-beta} for $a_i$ and the coefficient
identity $a_{ij}=\sum_r d^r_{ij}a_r$ established above, we have
\begin{align*}
 L(\lambda(a))
 &=\sum_i\pi(h^i)\beta(a_i)\\
 &=\sum_{i,j}\pi(h^i s(h^j))a_{ij}\\
 &=\sum_r\pi\!\left(m_K(\operatorname{id}_K\otimes s)\delta(h^r)\right)a_r\\
 &=\sum_r\varepsilon_K(h^r)a_r=a.
\end{align*}
This uses the antipode identity for $\delta$ and then the counit
identity for $\lambda$. For equivariance, use the equivariance of
$\pi$ and the invariance of $\beta(a)$ to compute
\begin{align*}
 \lambda^\infty(L(\phi\otimes a))
 &=\left(\sum_{(\phi)}\phi_{(2)}\otimes\pi(\phi_{(1)})\right)
   (1_K\otimes\beta(a))\\
 &=\sum_{(\phi)}\phi_{(2)}\otimes L(\phi_{(1)}\otimes a)
   =(\operatorname{id}_K\otimes L)\Gamma(\phi\otimes a).
\end{align*}
By linearity and continuity this holds on $K\otimes A$. Set $\psi=L\circ T_\lambda$. The map $\psi$ is a unital
$*$-homomorphism, and $\psi(1\boxtimes a)=L(\lambda(a))=a$.
Moreover,
\[
 \lambda^\infty\psi
 =(\operatorname{id}_K\otimes L)\Gamma T_\lambda
 =(\operatorname{id}_K\otimes (L T_\lambda))\Omega
 =(\operatorname{id}_K\otimes\psi)\Omega.
\]
This is the required equivariant sequential splitting. In particular,
for the ordered product $\phi\boxtimes a=(\phi\boxtimes1)(1\boxtimes a)$,
\[
 \psi(\phi\boxtimes a)=L((\phi\otimes1)\lambda(a))=\pi(\phi)a.
\]
This last formula is a consequence of the construction, not an
assumption about multiplication in the braided product.

\smallskip
\noindent\emph{From spatial Rokhlin property to Hopf Rokhlin property.}
No saturation is assumed in this direction.
Let $\psi:D\to A^\infty$ be an equivariant splitting of $j_A$.
Set $L=\psi\circ T_\lambda^{-1}$. Then $L$ is a $*$-homomorphism,
$L(\lambda(a))=\psi(1\boxtimes a)=a$, and
\begin{align*}
 \lambda^\infty L
 &=(\operatorname{id}_K\otimes\psi)\Omega T_\lambda^{-1}\\
 &=(\operatorname{id}_K\otimes(\psi T_\lambda^{-1}))\Gamma
   =(\operatorname{id}_K\otimes L)\Gamma.
\end{align*}
Since $\lambda(1_A)=1_K\otimes1_A$, the identity $L\lambda(1_A)=1_A$
shows that $L$ is unital. Define $\pi(\phi)=L(\phi\otimes1_A)$.
This is a unital $*$-homomorphism, as it is the composition of
$\phi\mapsto\phi\otimes1_A$ with $L$. Its equivariance follows from
\[
 \lambda^\infty(\pi(\phi))
 =(\operatorname{id}_K\otimes L)(\delta(\phi)\otimes1_A)
 =\sum_{(\phi)}\phi_{(2)}\otimes\pi(\phi_{(1)}).
\]
Evaluating the first leg at $h\in H$ gives
$h\cdot\pi(\phi)=\sum_{(\phi)}\phi_{(2)}(h)\pi(\phi_{(1)})$.

Finally, $e_K$ is a projection, so $p=\pi(e_K)$ is a projection.
To see that it commutes with $A$, use the distinguished projection
identities $e_K\chi=\chi e_K=\varepsilon_K(\chi)e_K$ for $\chi\in K$.
With $\lambda(a)=\sum_i h^i\otimes a_i$, both products
$(e_K\otimes1)\lambda(a)$ and $\lambda(a)(e_K\otimes1)$ are equal to
$e_K\otimes\sum_i\varepsilon_K(h^i)a_i=e_K\otimes a$ by the counit
identity. Applying $L$ and using $L\lambda(a)=a$, we obtain
\[
 \pi(e_K)a=L((e_K\otimes1)\lambda(a))
          =L(\lambda(a)(e_K\otimes1))=a\pi(e_K).
\]
Thus $\pi(e_K)\in A^\infty\cap A'$ and the two conditions in
Definition~\ref{def:RP-Hopf} hold.
\end{proof}

\subsection{Approximate representability}

We first define approximate representability for Hopf algebra actions,
and then recall the corresponding notion for conditional expectations.

\begin{definition}\label{AR-action}\cite[Definition~4.3]{KT-2015}
Let a finite-dimensional Hopf $\Cst$-algebra $H$ act on a unital
$\Cst$-algebra $A$, and let
\[
 \rho:A\longrightarrow A\otimes K,\qquad
 \rho(a)=\sum_i(h_i\cdot a)\otimes h^i,
\]
be the associated coaction from Proposition~\ref{thm: act-coact}.
The action is \emph{approximately representable} if there exists a
unitary $u\in A^\infty\otimes K$ satisfying the following conditions:
    \begin{enumerate}
        \item $\rho(a) = (Ad(u) \circ \rho_{K}^A) (a)$ for $a \in A$,
        \item $1 = (u \otimes 1) (\rho_{K}^{A^{\infty}} \otimes \operatorname{id})(u) (\operatorname{id} \otimes \Delta)(u^*)$,
        \item $1 = (\rho^{\infty} \otimes \operatorname{id})(u)(u \otimes 1)(\operatorname{id} \otimes \Delta)(u^*)$,
    \end{enumerate}
    where $\rho_{K}^{A^{\infty}}$ is the trivial coaction ($a \mapsto a \otimes 1$).
\end{definition}

\begin{remark}\label{rem:ar-leg-form}
Let
$u=\sum_i a_i\otimes h^i\in A^\infty\otimes K$ be a unitary.
The tensor-leg notation used below means
\[
 u_{12}=\sum_i a_i\otimes h^i\otimes1_K,
 \qquad
 u_{13}=\sum_i a_i\otimes1_K\otimes h^i
\]
in $A^\infty\otimes K\otimes K$.
Condition (1) of Definition~\ref{AR-action} is the implementation
identity
\[
 \rho(a)=u(a\otimes1_K)u^*\qquad(a\in A).
\]
For the remaining conditions, the trivial coaction gives
\[
 (\rho_K^{A^\infty}\otimes\operatorname{id})(u)
   =\sum_i a_i\otimes1_K\otimes h^i=u_{13}.
\]
$(\operatorname{id}\otimes\Delta_K)(u)$ is unitary because $\operatorname{id}\otimes\Delta_K$ is a
unital $*$-homomorphism. Thus condition (2) is
$1=u_{12}u_{13}(\operatorname{id}\otimes\Delta_K)(u)^*$. Multiplying on the right by $(\operatorname{id}\otimes\Delta_K)(u)$ shows that
it is equivalent to
\begin{equation}\label{eq:corep}
 (\operatorname{id}\otimes\Delta_K)(u)=u_{12}u_{13}.
\end{equation}
Likewise, condition (3) is
\[
 1=(\rho^\infty\otimes\operatorname{id})(u)u_{12}(\operatorname{id}\otimes\Delta_K)(u)^*.
\]
Multiplying on the right first by $(\operatorname{id}\otimes\Delta_K)(u)$ and then by $u_{12}^*$ gives
\[
 (\rho^\infty\otimes\operatorname{id})(u)= (\operatorname{id}\otimes\Delta_K)(u)u_{12}^*.
\]
Using \eqref{eq:corep}, this becomes
\begin{equation}\label{eq:equiv}
 (\rho^\infty\otimes\operatorname{id})(u)=u_{12}u_{13}u_{12}^*.
\end{equation}
Each multiplication is reversible. Hence, for a unitary $u$,
Definition~\ref{AR-action} is equivalent to the implementation
identity together with \eqref{eq:corep} and \eqref{eq:equiv}.
Here $\Delta_K$ is the original coproduct; no opposite coproduct
is involved in this reformulation.
\end{remark}

Following Lee and Osaka \cite[Definition~3.4]{LO2019}, we use the
following definition of approximate representability for conditional
expectations.

\begin{definition}\label{AR-incl}
A conditional expectation $E : A  \to  P$ of index-finite type is said to be \emph{approximately representable} if there exist a projection $ e \in P^{\infty} \cap P'$ and
a finite set of elements $\{u_i\} \subset A$ such that
\begin{enumerate}
    \item  $exe = E(x)e $ for every $x \in A$,
    \item $\sum_i u_i e u_i^* = 1$,
    \item the map $x \mapsto xe$ is injective for $x \in P$.
\end{enumerate}
\end{definition}

For a finite group $G$, an action of the Hopf $\Cst$-algebra
$(C_r^*(G),\Delta)$ on $A$ corresponds to a group action $\alpha$ of $G$
on $A$. In this case, Definition~\ref{AR-action} is the condition suggested by
Izumi for general finite groups in \cite[Remark~3.7]{Izumi2004I}.
When $G$ is abelian and both $A$ and $A\rtimes_\alpha G$ are
simple, the equivalence in Theorem~\ref{thm:main} is the one established in
\cite[Proposition~3.17]{LO2019}. The corresponding equivalence for
the Rokhlin property holds for actions of arbitrary finite groups on
unital simple $\Cst$-algebras \cite[Proposition~3.12]{LO2019}.

\begin{theorem}\label{thm:main}
Let $H\curvearrowright A$ be a saturated action of a finite-dimensional
Hopf $\Cst$-algebra on a unital $\Cst$-algebra. Let
$E\colon A\to P$ be the canonical conditional expectation. If $P$ is
simple, then the action is approximately representable if and only if
$E$ is approximately representable. The implication from approximate
representability of $E$ to approximate representability of the action
holds without any simplicity assumption.
\end{theorem}

\begin{proof}
\allowdisplaybreaks[2]
\noindent\emph{From the action to the expectation.}
We use the identification $H\cong(H^\circ)^\circ = K^\circ$. For convenience, write
$\varepsilon$, $S$, and $\Delta$ for the structure maps on $K$.
Choose an implementing unitary $u\in A^\infty\otimes K$ as in
Definition~\ref{AR-action}. By Remark~\ref{rem:ar-leg-form},
$\rho(a)=u(a\otimes1)u^*$ for $a\in A$, and $u$ satisfies
\eqref{eq:corep} and \eqref{eq:equiv}.

Write $u = \sum_i a_i \otimes h^i $ for some $a_i \in A^{\infty}$, and set $v := (\operatorname{id} \otimes \varepsilon)(u)$. Applying $(\operatorname{id} \otimes \varepsilon \otimes \operatorname{id})$ to \eqref{eq:corep} and using $(\operatorname{id} \otimes \varepsilon)\Delta = \operatorname{id}$, we obtain $u = (v \otimes 1)u$. Applying $(\operatorname{id} \otimes m(S \otimes \operatorname{id}))$ to \eqref{eq:corep}, we obtain $(\operatorname{id} \otimes \varepsilon)(u) = (\sum_i a_i \otimes S(h^i)) u$, which implies $1 = (\operatorname{id} \otimes S)(u)(u)$. Therefore, the two equations above give
\begin{equation*}
   (\operatorname{id} \otimes \varepsilon)(u) = 1 \text{ and } u^* = (\operatorname{id} \otimes S)(u).
\end{equation*} 
Define the map
\[
  \pi : H \longrightarrow A^{\infty}, \qquad
  \pi(h) := (\operatorname{id} \otimes h)(u).
\]

The multiplicativity of $\pi$ follows from
\eqref{eq:corep}, as
\[
  \pi(h)\pi(k)
  = (\operatorname{id} \otimes h \otimes k)(u_{12}u_{13})
  = (\operatorname{id} \otimes h \otimes k)(\operatorname{id} \otimes \Delta)(u)
  = \pi(hk).
\]
Since
$1_{H} = \varepsilon$ and $(\operatorname{id} \otimes \varepsilon)(u) = 1$, $\pi(1_{H}) = \pi(\varepsilon) = (\operatorname{id} \otimes \varepsilon) u = 1_{A^{\infty}}$. For $y\in K$, write $\omega_h(y) := \overline{h(y^{*})}$, so that $h^{*}(y) = \overline{h(S(y)^{*})} = \omega_h(S(y))$.
Then
\[
  \pi(h)^{*} = \sum_i a_i^* \overline{h(h^i)} 
  = (\operatorname{id} \otimes \omega_h)(u^{*})
  = (\operatorname{id} \otimes \omega_h \circ S)(u)
  = (\operatorname{id} \otimes h^{*})(u)
  = \pi(h^{*}).
\]

\textbf{Claim 1.} \emph{For all $h \in H$ and $a \in A$,}
\begin{equation}\label{eq:pi ba}
  \pi(h)\,a \;=\; \sum (h_{(1)} \cdot a)\,\pi(h_{(2)}).
\end{equation}

\smallskip
To obtain \eqref{eq:pi ba}, slice $\rho(a)u=u(a\otimes 1)$ by $h$.
On the right, the slice is $\pi(h)a$. On the left, multiplication
in $K$ is dual to the coproduct $\Delta_H$, so slicing the two
factors by $h_{(1)}$ and $h_{(2)}$ gives the sum in \eqref{eq:pi ba}.
Let $e=\pi(e_H)$. This is a projection, since $\pi$ is a
$*$-homomorphism. We will show that it is nonzero in Claim~4.

\textbf{Claim 2.} \emph{$e\in P^\infty\cap P'$.}
For $x\in P$, we have $h\cdot x=\varepsilon_H(h)x$. Applying this
identity in \eqref{eq:pi ba} and using the counit identity gives
\begin{align*}
 \pi(h)x
 &=\sum_{(h)}(h_{(1)}\cdot x)\pi(h_{(2)})
   =x\sum_{(h)}\varepsilon_H(h_{(1)})\pi(h_{(2)})\\
 &=x\,\pi\left(\sum_{(h)}\varepsilon_H(h_{(1)})h_{(2)}\right)
   =x\pi(h).
\end{align*}
In particular, $e\in P'$. Because $e_H$ is central in $H$, $e$ commutes
with each $\pi(h_i)$. Thus
\begin{equation*}
 u(e\otimes 1)=\sum_i\pi(h_i)e\otimes h^i
             =\sum_i e\pi(h_i)\otimes h^i=(e\otimes 1)u.
\end{equation*}
Slicing the third leg of \eqref{eq:equiv} by $e_H$ therefore gives
\[
 \rho^\infty(e)=u(e\otimes 1)u^*=e\otimes 1.
\]
In particular, $E^\infty(e)=e$. Since $E^\infty$ takes values in
$P^\infty$, this gives $e\in P^\infty\cap P'$.

\textbf{Claim 3.} \emph{$exe=E(x)e$ for every $x\in A$.}
Apply \eqref{eq:pi ba} with $h=e_H$ and use the defining relation
$(e_H)_{(2)}e_H=\varepsilon_H((e_H)_{(2)})e_H$ in the sum. Then
\begin{align*}
 exe&=\pi(e_H)x\pi(e_H)\\
 &=\sum_{(e_H)}((e_H)_{(1)}\cdot x)\,
                   \pi((e_H)_{(2)}e_H)\\
 &=\sum_{(e_H)}\varepsilon_H((e_H)_{(2)})
                   ((e_H)_{(1)}\cdot x)e\\
 &=(e_H\cdot x)e=E(x)e,
\end{align*}
where the penultimate equality is the counit identity.

\textbf{Claim 4.}  \emph{There is a finite set $\{u_{i}\} \subset A$
with $\sum_{i} u_{i}eu_{i}^{*} = 1$.}

Using $\pi$, define the linear map $\Phi : A \rtimes H \to A^{\infty}$
by $\Phi(a \otimes h) := a \pi(h)$. Then
$\Phi|_{A} = \operatorname{id}_{A}$ and $\Phi(1 \otimes e_H) = \pi(e_H)$. Note that $\Phi$ is a unital $*$-homomorphism because
\begin{align*}
    \Phi((a \otimes h)( b \otimes k)) & = \Phi \left(\sum a(h_{(1)} \cdot b) \otimes h_{(2)} k\right) = \sum a(h_{(1)} \cdot b) \pi (h_{(2)} k) \\
    & = \sum a(h_{(1)} \cdot b ) \pi (h_{(2)}) \pi(k) = a \left(\sum (h_{(1)} \cdot b) \pi(h_{(2)})\right) \pi(k)\\
    & = a \pi(h) b \pi(k)  = \Phi(a \otimes h) \Phi(b \otimes k) 
\end{align*}

\begin{align*}
    \text{ and } \Phi (a \otimes h)^* & = (a \pi(h))^*  = \pi(h)^* a^*  = \pi(h^*) a^*  = \sum (h_{(1)}^* \cdot a^*) \pi(h_{(2)}^*)\\
    & = \Phi \left(\sum (h_{(1)}^* \cdot a^*) \otimes h_{(2)}^*\right) = \Phi((a \otimes h)^*).
\end{align*}

Since $\pi$ is unital, $\Phi$ is also unital. $E$ is of index-finite type, so $E$ admits a quasi-basis $\{(u_{i},u_{i}^{*})\}_{i=1}^{n}$ with $u_{i} \in A$, and
$P \subset A \subset A \rtimes H$ is a basic construction with Jones
projection $e_{1} = 1 \otimes e_H$ (see \cite{SP1994}) and $\sum_{i} u_{i}e_{1}u_{i}^{*} = 1$ in $A \rtimes H$. Applying the unital
$*$-homomorphism $\Phi$, which fixes each $u_{i} \in A$ and carries
$e_{1}$ to $e$, gives
\[
  \sum_{i} u_{i}\,e\,u_{i}^{*} = \Phi(1) = 1 .
\]

\medskip
The identity $\sum_i u_i e u_i^*=1$ proves that $e\ne0$.
It remains to verify condition (3) of Definition~\ref{AR-incl}. The map
$P\ni x\mapsto xe$ is a $*$-homomorphism because $e$ is a projection
commuting with $P$. Its kernel is a closed two-sided ideal of $P$.
This ideal is proper, since $1e=e\ne0$, and it is therefore zero by
simplicity of $P$. This is the only use of simplicity in this direction.
Hence $E$ is approximately representable.

\smallskip
\noindent\emph{From the expectation to the action.}
No simplicity assumption is needed in this direction.

\smallskip
\noindent\emph{Preparing the canonical unitary $W$.} Define
\[
 W=\sum_i h_i\otimes h^i\in H\otimes K.
\]
For $x\in H$, write $\operatorname{ev}_x(\phi)=\phi(x)$. The dual
involution and the dual-basis expansion give
\[
 (\operatorname{id}\otimes\operatorname{ev}_x)(W^*)
 =\left(\sum_i h^i(S_H(x)^*)h_i\right)^*=S_H(x).
\]
Thus $W^*=(S_H\otimes\operatorname{id})(W)$. Since multiplication
in $K$ is dual to $\Delta_H$, the antipode identities give
\begin{align*}
 (\operatorname{id}\otimes\operatorname{ev}_x)(W^*W)
 &=\sum_{(x)}S_H(x_{(1)})x_{(2)}=\varepsilon_H(x)1_H,\\
 (\operatorname{id}\otimes\operatorname{ev}_x)(WW^*)
 &=\sum_{(x)}x_{(1)}S_H(x_{(2)})=\varepsilon_H(x)1_H.
\end{align*}
These are the slices of $1_H\otimes1_K$. Since the evaluations
separate points, $W$ is unitary.

\smallskip
\noindent\emph{Constructing the map $\Phi$.}
Let $e_{1} = 1 \otimes e_H$ be the Jones projection as before. Since $E $ is approximately representable, by Definition~\ref{AR-incl}, there exist a projection $e \in P^{\infty} \cap P'$ and a finite set $\{u_i:  i \in I\}$
satisfying the conditions in Definition~\ref{AR-incl}. We need to construct a suitable unitary $u$ that satisfies the conditions in Definition~\ref{AR-action}. Since the action is saturated, $A \rtimes H = \operatorname{span}\{x e_1 y : x, y \in A \}$; the algebraic span suffices, since it is an ideal containing $1$. Define the map
\[
\Phi : A \rtimes H \to A^{\infty} \text{ such that }
\sum_{l} x_{l}e_{1}y_{l} \mapsto \sum_{l} x_{l}\,e\,y_{l}.
\]
This is the map in \cite[Proposition~3.6]{LO2019} transported to the
crossed product. For clarity, we verify independence of the chosen
expression. Suppose $\sum_l x_l e_1 y_l=0$. Acting on $z\in A$ in the
basic-construction representation gives $\sum_l x_l E(y_lz)=0$, since
$e_1$ acts there as $E$. Put $t=\sum_l x_l e y_l$. The covariance
relation $eae=E(a)e$ implies
\[
 tze=\sum_l x_l e y_lze=\sum_l x_l E(y_lz)e=0\qquad(z\in A).
\]
Using the family in Definition~\ref{AR-incl}, we obtain
$t=t\sum_i u_i e u_i^*=\sum_i(tu_i e)u_i^*=0$. Thus $\Phi$ is
well defined and linear. If $\{(v_j,v_j^*)\}$ is a quasi-basis for
$E$, then $\sum_j v_j e_1v_j^*=1$. To see that $\Phi$ fixes $A$,
first note that its range contains $1=\sum_i u_i e u_i^*$. The
multiplicativity proved below therefore implies $\Phi(1)=1$.
Then, writing $a=\sum_j av_j e_1v_j^*$ gives
$\Phi(a)=a\sum_j v_j e v_j^*=a\Phi(1)=a$. For products, $(xe_{1}y)(ze_{1}w) = x\,E(yz)\,e_{1}w$ in $A \rtimes H$ because
$e_{1}ze_{1} = E(z)e_{1}$, and we have $(xey)(zew) = xE(yz)ew$ in
$A^{\infty}$ by Definition~\ref{AR-incl}, so $\Phi$ is multiplicative. It preserves $*$ because
$(xe_{1}y)^{*} = y^{*}e_{1}x^{*}$ and $(xey)^{*} = y^{*}ex^{*}$. 

\smallskip
\noindent\emph{Constructing the implementing unitary $u$.}
Set
\[
u:= \sum_i \Phi(h_i) \otimes h^i = (\Phi \otimes \operatorname{id})(W).
\] 
Since $W$ is unitary, so is $u$. Here $\Phi(h)$ means
$\Phi(1\otimes h)$ under the canonical embedding
$H\hookrightarrow A\rtimes H$. The dual-basis expansion gives
$(\operatorname{id}\otimes h)(u)=\Phi(h)$ for every $h\in H$.

\smallskip
\noindent\emph{Verifying the defining conditions for $u$.}
We check the three conditions in Definition~\ref{AR-action} in order. 

\smallskip
\noindent\emph{Condition (1): implementation.}
Write $h\in H$ also for its evaluation functional on $K$.
Multiplication in $K$ is dual to $\Delta_H$, so
$(\phi\chi)(h)=\sum_{(h)}\phi(h_{(1)})\chi(h_{(2)})$.
Using this identity and the dual-basis expansion, followed by the
crossed-product relation and $\Phi|_A=\operatorname{id}_A$, gives
\begin{align}
 (\operatorname{id}\otimes h)(\rho(a)u)
 &=\sum_{(h)}(h_{(1)}\cdot a)\Phi(h_{(2)})\notag\\
 &=\Phi\left(\sum_{(h)}(h_{(1)}\cdot a)\otimes h_{(2)}\right)\notag\\
 &=\Phi\bigl((1\otimes h)(a\otimes1)\bigr)
   =\Phi(h)a\notag\\
 &=(\operatorname{id}\otimes h)(u(a\otimes1)).\label{AR-prop1}
\end{align}
The evaluations separate points of $K$, so
$\rho(a)u=u(a\otimes1)$.

\smallskip
\noindent\emph{Condition (2): corepresentation identity.}
Since $\Delta_K$ is dual to multiplication in $H$, slicing at
$h,k\in H$ and using multiplicativity of $\Phi$ gives
\begin{align}
 (\operatorname{id}\otimes h\otimes k)
       ((\operatorname{id}\otimes\Delta_K)(u))
 &=\Phi(hk)=\Phi(h)\Phi(k)\notag\\
 &=(\operatorname{id}\otimes h\otimes k)(u_{12}u_{13}).\label{AR-prop2}
\end{align}
These slices separate points, giving the corepresentation identity.

\smallskip
\noindent\emph{Condition (3): equivariance of $u$.}
By condition (2) and unitarity, condition (3) is equivalent to
$(\rho^\infty\otimes\operatorname{id})(u)=u_{12}u_{13}u_{12}^*$.
We first show that $\Phi$ intertwines the action on $A^\infty$ with
the extended action on $A\rtimes H$ (see the discussion after
Proposition~3.4 in \cite{SP1994}):
\[
 h\cdot z=\sum_{(h)}(1\otimes h_{(1)})z(1\otimes S_H(h_{(2)})).
\]
Indeed, $e_1=1\otimes e_H$ transforms by the counit. The Haar identity
$h_{(1)}e_H=\varepsilon_H(h_{(1)})e_H$, the counit identity, and
$e_H S_H(h)=\varepsilon_H(h)e_H$ give
\begin{align*}
 h\cdot e_1
 &=1\otimes\sum_{(h)}h_{(1)}e_H S_H(h_{(2)})\\
 &=1\otimes e_H S_H\left(\sum_{(h)}\varepsilon_H(h_{(1)})h_{(2)}\right)\\
 &=1\otimes e_H S_H(h)=\varepsilon_H(h)e_1.
\end{align*}
The projection $e$ also transforms by the counit, since $e\in P^\infty$
and $P=A^H$. To apply the module-algebra rule
$h\cdot(xy)=\sum_{(h)}(h_{(1)}\cdot x)(h_{(2)}\cdot y)$ twice,
write the iterated coproduct as
$(\operatorname{id}_H\otimes\Delta_H)\Delta_H(h)=
\sum_{(h)}h_{(1)}\otimes h_{(2)}\otimes h_{(3)}$.
Using $\Phi(xe_1y)=xey$ and the counit in the middle leg gives
\begin{align*}
 \Phi(h\cdot(xe_1y))
 &=\sum_{(h)}\varepsilon_H(h_{(2)})
     \Phi\bigl((h_{(1)}\cdot x)e_1(h_{(3)}\cdot y)\bigr)\\
 &=\sum_{(h)}\varepsilon_H(h_{(2)})
     (h_{(1)}\cdot x)e(h_{(3)}\cdot y)\\
 &=h\cdot(xey)=h\cdot\Phi(xe_1y).
\end{align*}
Since the elements $xe_1y$ span $A\rtimes H$ by saturation, this proves
that $\Phi$ is equivariant.

To compare slices, use
$u^*=\sum_l\Phi(S_H(h_l))\otimes h^l$.
Multiplication in $K$ is dual to $\Delta_H$, as in
\eqref{AR-prop1}. Using this fact, the extended adjoint action,
and the equivariance of $\Phi$ just proved, we obtain
\begin{align*}
 &(\operatorname{id}\otimes h\otimes k)(u_{12}u_{13}u_{12}^*)\\
 &\quad=\sum_{(h)}\Phi(h_{(1)})\Phi(k)\Phi(S_H(h_{(2)}))\\
 &\quad=\Phi\left(\sum_{(h)}h_{(1)}kS_H(h_{(2)})\right)\\
 &\quad=\Phi(h\cdot(1\otimes k))=h\cdot\Phi(k)\\
 &\quad=(\operatorname{id}\otimes h\otimes k)
             ((\rho^\infty\otimes\operatorname{id})(u)).
\end{align*}
Since $h,k$ are arbitrary, we obtain
\begin{equation}\label{AR-prop3}
 u_{12}u_{13}u_{12}^*=(\rho^\infty\otimes\operatorname{id})(u).
\end{equation}
Thus, by \eqref{AR-prop1}, \eqref{AR-prop2}, and \eqref{AR-prop3}, the action of $H$ on $A$ is approximately representable.

\end{proof}

To compare approximate representability with its spatial formulation,
we first make the tensor-leg convention explicit. As before $C(\mathbb G)=K$, and
$\delta=\Delta_{\mathbb G}=\Delta_K^{\mathrm{cop}}$ as above, and put
$\lambda=\Sigma\rho$. The coproduct in Definition~\ref{AR-action}
is the original $\Delta_K$; the coproduct for the action $\lambda$ is
$\delta$. Thus the coproduct being flipped belongs to the coaction
algebra $K$, not to the acting algebra $H$.

\begin{lemma}\label{lem:ar-leg-convention}
Let $H$ act on a unital $\Cst$-algebra $A$, and use the notation above.
A unitary $u\in A^\infty\otimes K$ implements approximate
representability as in Definition~\ref{AR-action} if and only if
$V=\Sigma(u^*)\in K\otimes A^\infty$ satisfies
\begin{equation}\label{eq:ar-bsv-relations}
 \begin{aligned}
  V^*(1_K\otimes a)V&=\lambda(a)\quad(a\in A),\\
  (\delta\otimes\operatorname{id})(V)&=V_{13}V_{23},\\
  (\operatorname{id}\otimes\lambda^\infty)(V)
    &=V_{23}^*V_{13}V_{23}.
 \end{aligned}
\end{equation}
The correspondence is reversible, with inverse
$u=\Sigma(V^*)$.
\end{lemma}

\begin{proof}
By Remark~\ref{rem:ar-leg-form}, approximate representability is
equivalent to implementation of $\rho$ together with
\eqref{eq:corep} and \eqref{eq:equiv}.

Besides the two-factor flip $\Sigma$, use the $*$-isomorphism
\[
 \mathcal T:A^\infty\otimes K\otimes K
 \longrightarrow K\otimes K\otimes A^\infty,\qquad
 \mathcal T(a\otimes\phi\otimes\chi)=\chi\otimes\phi\otimes a.
\]
Thus $\mathcal T$ preserves multiplication order; taking adjoints,
rather than permuting tensor factors, reverses the order of a
product. For $V=\Sigma(u^*)$, inspection of the tensor legs gives
\[
 \mathcal T(u_{12}^*)=V_{23},\qquad
 \mathcal T(u_{13}^*)=V_{13}.
\]
Moreover, for every $z\in A^\infty\otimes K$,
\[
 \begin{aligned}
 \mathcal T((\operatorname{id}\otimes\Delta_K)(z))
   &=(\delta\otimes\operatorname{id})(\Sigma(z)),\\
 \mathcal T((\rho^\infty\otimes\operatorname{id})(z))
   &=(\operatorname{id}\otimes\lambda^\infty)(\Sigma(z)).
 \end{aligned}
\]
Indeed, the first identity sends $a\otimes\phi$ to
$\sum_{(\phi)}\phi_{(2)}\otimes\phi_{(1)}\otimes a$, using
$\delta=\Sigma\Delta_K$; the second follows from
$\lambda^\infty=\Sigma\rho^\infty$. Linearity gives the identities
for arbitrary $z$.

Suppose first that $u$ implements approximate representability.
Since $\Sigma$ is a $*$-isomorphism, $V$ is unitary, and flipping the
implementation identity gives
\[
 \lambda(a)=\Sigma\bigl(u(a\otimes1_K)u^*\bigr)
           =V^*(1_K\otimes a)V.
\]
Taking adjoints in \eqref{eq:corep} and \eqref{eq:equiv}, and then
applying $\mathcal T$, gives
\[
 \begin{aligned}
 (\delta\otimes\operatorname{id})(V)
   &=\mathcal T(u_{13}^*u_{12}^*)=V_{13}V_{23},\\
 (\operatorname{id}\otimes\lambda^\infty)(V)
   &=\mathcal T(u_{12}u_{13}^*u_{12}^*)
     =V_{23}^*V_{13}V_{23}.
 \end{aligned}
\]
These are precisely \eqref{eq:ar-bsv-relations}.

Conversely, suppose that $V$ is unitary and satisfies
\eqref{eq:ar-bsv-relations}, and set $u=\Sigma(V^*)$.
Then $u$ is unitary. Flipping the first relation gives
$\rho(a)=u(a\otimes1_K)u^*$. Applying $\mathcal T^{-1}$ to the
other two relations, using the identities above, gives
\[
 \begin{aligned}
 (\operatorname{id}\otimes\Delta_K)(u^*)
   &=u_{13}^*u_{12}^*,\\
 (\rho^\infty\otimes\operatorname{id})(u^*)
   &=u_{12}u_{13}^*u_{12}^*.
 \end{aligned}
\]
Taking adjoints recovers \eqref{eq:corep} and \eqref{eq:equiv}.
Remark~\ref{rem:ar-leg-form} therefore implies that $u$ implements
approximate representability.
\end{proof}

We recall the spatial definition before applying the lemma. Regard the
finite quantum group $\mathbb G$ also as a discrete quantum group, and
put $\mathcal C=\mathbb G\ltimes_{\lambda,\mathrm r}A$.
Let $\iota_{\mathbb G}:C_r^*(\mathbb G)\to\mathcal C$ be the canonical
embedding and put
$W_\lambda=(\operatorname{id}_K\otimes\iota_{\mathbb G})(W_{\mathbb G})$,
where $W_{\mathbb G}$ is the multiplicative unitary of $\mathbb G$.
By \cite[Definition~4.7]{BSV2017}, spatial approximate representability means
that the canonical embedding $j_\lambda:A\to\mathcal C$, with
$\mathcal C$ carrying the inner coaction
\[
 \gamma_\lambda(x)=W_\lambda^*(1_K\otimes x)W_\lambda,
\]
is equivariantly sequentially split. In the present unital setting,
this means that there is a $*$-homomorphism $\psi:\mathcal C\to A^\infty$
such that
\[
 \psi(j_\lambda(a))=a,\qquad
 \lambda^\infty\psi=(\operatorname{id}_K\otimes\psi)\gamma_\lambda.
\]

\Needspace{9\baselineskip}
\begin{theorem}\label{thm:AR-SAR}
Let a finite-dimensional Hopf $\Cst$-algebra $H$ act on a separable
unital $\Cst$-algebra $A$. Let
$C(\mathbb G)=K$ with $\Delta_{\mathbb G}=\Delta_K^{\mathrm{cop}}$,
and let $\lambda=\Sigma\rho$ be the associated action of
$\mathbb G$ on $A$. Then $H\curvearrowright A$ is approximately
representable if and only if $\lambda$ is spatially approximately
representable in the sense of Barlak--Szab\'o--Voigt.
Neither saturation nor simplicity is required for this comparison.
\end{theorem}

\begin{proof}
The quantum group $\mathbb G$ is both compact and discrete, since
$C(\mathbb G)$ and $C_r^*(\mathbb G)$ are finite-dimensional.
It is exact because its reduced group $\Cst$-algebra is
finite-dimensional, and it is amenable because it is compact;
see \cite[Definitions~1.7,~1.12, and~2.7]{BSV2017} and the observation
following Definition~1.12 there. Thus both directions of
\cite[Proposition~4.8]{BSV2017} apply. Moreover, in the notation of that
proposition, $D_{\infty,A}=A^\infty$ because $A$ is unital, and
$M(K\otimes A^\infty)=K\otimes A^\infty$.

Suppose first that $H\curvearrowright A$ is approximately
representable, with implementing unitary $u$.
Lemma~\ref{lem:ar-leg-convention} gives a unitary $V=\Sigma(u^*)$
satisfying \eqref{eq:ar-bsv-relations}. The second relation says
precisely that $V$ is a unitary representation of $\mathbb G$ in the
sense of \cite[Definition~1.8]{BSV2017}. The first and third relations
are the implementation and equivariance requirements in
\cite[Proposition~4.8(b)]{BSV2017}. That proposition therefore gives the
equivariant sequential splitting in the spatial definition.

Conversely, if $\lambda$ is spatially approximately representable,
\cite[Proposition~4.8(a)]{BSV2017} provides such a unitary representation
$V\in K\otimes A^\infty$. It satisfies the first and third relations
in \eqref{eq:ar-bsv-relations}; its representation identity is the
second relation. By the converse direction of
Lemma~\ref{lem:ar-leg-convention}, $u=\Sigma(V^*)$ satisfies
Definition~\ref{AR-action}. Hence $H\curvearrowright A$ is
approximately representable.
\end{proof}

\subsection{Duality for saturated Hopf actions}
We combine the comparisons to prove the
dynamical duality from the inclusion alone. The key is that the
fixed-point expectation for the dual action is the dual conditional expectation of
the original fixed-point expectation. This identifies the mechanism
that exchanges the Rokhlin property and approximate representability.

\begin{theorem}\label{Rok-AR-duality}
Let $H\curvearrowright A$ be a saturated action of a finite-dimensional
Hopf $\Cst$-algebra $H$ on a unital simple $\Cst$-algebra $A$. Suppose
that $B=A\rtimes H$ is simple. Then
\begin{enumerate}
 \item The action of $H$ on $A$ has the Rokhlin property if and only if
       the dual action of $K$ on $B$ is approximately representable.
 \item The action of $H$ on $A$ is approximately representable if and
       only if the dual action of $K$ on $B$ has the Rokhlin property.
\end{enumerate}
\end{theorem}

\begin{proof}
\allowdisplaybreaks[2]
Let $N=\dim H$, and let $E\colon A\to P$ be the canonical
expectation. Saturation identifies the basic construction for $E$
with $B$, with Jones projection $e_1=1\otimes e_H$, and gives
$\operatorname{Index}(E)=N1_A$.
The relation $e_1ae_1=E(a)e_1$ and the spanning family $ae_1b$ show
that $e_1Be_1=Pe_1$. The map $p\mapsto pe_1$ is faithful on $P$:
in the basic-construction representation, $pe_1$ sends $1_A$ to $p$.
Therefore $P\cong e_1Be_1$ is simple because $B$ is simple.

\smallskip
\noindent\emph{The dual action and the dual expectation.}
The dual coaction on $B$ is given on elementary tensors by
\[
 \widehat\rho(a\otimes h)
  =\sum_{(h)}(a\otimes h_{(1)})\otimes h_{(2)},
\]
and the dual action of $K$ is obtained by slicing its last leg.
The distinguished projection $e_K$ is the Haar trace $\tau_H$ viewed
as an element of $K$. Hence the canonical expectation of the dual
action is
\begin{align*}
 F(a\otimes h)
 &=e_K\cdot(a\otimes h)\\
 &=\sum_{(h)}(a\otimes h_{(1)})\tau_H(h_{(2)})\\
 &=\tau_H(h)(a\otimes 1_H),
\end{align*}
where the last equality is Haar invariance. Its range is $A$, so $B^K=A$. This is the canonical expectation
$F:B\to A$ introduced in the preliminaries. By
Proposition~\ref{prop:crossed-product-basic-construction}, it is
$\widehat E$ under the basic-construction identification.

We verify saturation of the dual action. Choose a quasi-basis
$\{(v_i,v_i^*)\}$ for $E$ and set $s_i=\sqrt N\,v_i e_1$. For
$x=ae_1b$ the Jones relation and the formula for $F$ give
\[
 s_i^*x=\sqrt N\,E(v_i^*a)e_1b,
 \qquad F(s_i^*x)=N^{-1/2}E(v_i^*a)b.
\]
Since $e_1$ commutes with $P$, the quasi-basis identity for $E$ yields
\[
 \sum_i s_i F(s_i^*x)
 =\sum_i v_i E(v_i^*a)e_1 b=ae_1b=x.
\]
The same identity holds for all $x\in B$ by linearity and continuity.
Taking adjoints after replacing $x$ by $x^*$ gives the other
quasi-basis identity. Thus $\{(s_i,s_i^*)\}$ is a quasi-basis for $F$,
and
\[
 \operatorname{Index}(F)=\sum_i s_i s_i^*
   =N\sum_i v_i e_1 v_i^*=N1_B.
\]
Since $\dim K=N$, the dual action is saturated by
\cite[Theorem~3.3]{JP2008}.

\smallskip
\noindent\emph{Combining the comparisons.}
We use the duality theorem of Lee and Osaka
\cite[Theorem~3.16]{LO2019}: a conditional expectation of index-finite type has the
Rokhlin property if and only if its dual is approximately
representable, and is approximately representable if and only if
its dual has the Rokhlin property. Both equivalences apply to $E$
and $F=\widehat E$.

For (1), simplicity of $A$ and Theorem~\ref{thm:Hopf-incl} identify the
Rokhlin property for $H\curvearrowright A$ with that of $E$.
Lee--Osaka duality identifies this with approximate representability
of $F$. Since $B^K=A$ is simple, Theorem~\ref{thm:main} for the saturated
dual action identifies approximate representability of $F$ with
that of $K\curvearrowright B$.

For (2), simplicity of $P$ and Theorem~\ref{thm:main} identify
approximate representability for $H\curvearrowright A$ with that
of $E$. By Lee--Osaka duality this is equivalent to the Rokhlin
property for $F$, which is equivalent to the Rokhlin property for
$K\curvearrowright B$ by Theorem~\ref{thm:Hopf-incl}, since $B$ is simple.
\end{proof}

\begin{remark}[The spatial approach]
If $A$ is separable, the two action-level equivalences in
Theorem~\ref{Rok-AR-duality} can also be deduced from Barlak--Szab\'o--Voigt.
We use the finite quantum group $\mathbb G$ with
$(K,\Delta_K^{\mathrm{cop}})$ and the left coaction
$\lambda=\Sigma\rho$ as in Theorem~\ref{thm:RP-SRP}.

To match the dual-action conventions, represent $K$ by left
multiplication $M$ on $L^2(K,\tau_K)$ and let
$(R_h\phi)(x)=\phi(xh)$ for $h,x\in H$ and $\phi\in K$. The
multiplicative unitary in \cite[Remarks~1.3--1.4]{BSV2017} satisfies
$W^*=\sum_i M_{h^i}\otimes R_{h_i}$: indeed,
$W^*(\phi\otimes\psi)=\Delta_K^{\mathrm{cop}}(\psi)(\phi\otimes1)$.
Its slices identify $C_r^*(\mathbb G)$ with $R(H)$. The coproduct on
$C(\check{\mathbb G})$ satisfies
$(\operatorname{id}\otimes\Delta_{\check{\mathbb G}})(W^*)=W_{13}^*W_{12}^*$
by the pentagon identity. Slicing the first leg by $h$ gives
$\Delta_{\check{\mathbb G}}(R_h)=\sum_{(h)}R_{h_{(2)}}\otimes R_{h_{(1)}}$.
The covariance relation
\[
 (R_h\otimes1)\lambda(a)
 =\sum_{(h)}\lambda(h_{(1)}\cdot a)(R_{h_{(2)}}\otimes1)
\]
follows from $R_hM_\phi=\sum_{(h)}M_{R_{h_{(1)}}\phi}R_{h_{(2)}}$,
which is the dual product formula evaluated at $xh$. Thus the map
$a\otimes h\mapsto\lambda(a)(R_h\otimes1)$ gives the crossed-product
map into the algebra of \cite[Definition~1.18]{BSV2017}. It is onto
by that definition and intertwines the dual coactions on the
original algebra and on $R(H)$ with $\Sigma\widehat\rho$ above.
It is faithful because it intertwines the dual Haar expectations,
is faithful on $A$, and $F$ is faithful. In particular the spatial
dual action is precisely the flipped dual coaction in the convention
of Theorem~\ref{thm:RP-SRP}.

The AR comparison, including the change of tensor-leg convention,
is proved in Lemma~\ref{lem:ar-leg-convention} and
Theorem~\ref{thm:AR-SAR}. Together with Theorem~\ref{thm:RP-SRP}
for the original and the dual actions, this identifies the two
properties with their spatial counterparts. The two exchanges now
follow from \cite[Theorem~4.12]{BSV2017}. The preceding proof instead
explains the exchanges through $E$ and its dual conditional expectation.
\end{remark}

\section{Examples, permanence, and obstructions}\label{sec:applications}

The first example applies our comparisons to the construction of
\cite[Example~5.9]{LO2026WT}. We recall the finite-group dictionary
in this context and refer to that preprint for the original
construction and its non-group type conclusion.

\begin{example}\label{ex:non-group-type}
Let $G$ be a finite nonabelian group with identity $1_G$, and put
$m=|G|$. The dual Hopf $\Cst$-algebras are
\[
 H=C^*(G)=\operatorname{span}\{s_g:g\in G\},\qquad
 K=C(G)=\operatorname{span}\{\delta_g:g\in G\},
\]
with pairing $\langle\delta_g,s_h\rangle=\delta_{g,h}$ and structure
maps
\begin{align*}
 \Delta_H(s_g)&=s_g\otimes s_g, &
 \Delta_K(\delta_k)&=\sum_{gh=k}\delta_g\otimes\delta_h,\\
 \varepsilon_H(s_g)&=1, &
 \varepsilon_K(\delta_g)&=\delta_{g,1_G},\\
 S_H(s_g)&=s_{g^{-1}}, & S_K(\delta_g)&=\delta_{g^{-1}}.
\end{align*}
An action $\alpha:G\curvearrowright A$ gives the $H$-action
$s_g\cdot a=\alpha_g(a)$ and the coaction
\[
 \rho_\alpha(a)=\sum_{g\in G}\alpha_g(a)\otimes\delta_g,
 \qquad
 \alpha_g(a)=(\operatorname{id}_A\otimes\operatorname{ev}_g)
                 \rho_\alpha(a),
\]
where $\operatorname{ev}_g(f)=f(g)$.

Now let $\lambda_g$ be the left regular representation of $G$ and set
\[
 A=\bigotimes_{n=1}^\infty M_m(\mathbb C),\qquad
 \alpha_g=\bigotimes_{n=1}^\infty\operatorname{Ad}(\lambda_g).
\]
This action is approximately representable and pointwise outer by
\cite[Example~5.9]{LO2026WT}.
Its implementing representation is
$v_g=[(\lambda_g^{\otimes n}\otimes1)_n]\in A^\infty$.
The representation and conjugacy-covariance relations from that
example show that $U=\sum_gv_g\otimes\delta_g$ satisfies
Definition~\ref{AR-action}. Hence the associated $H$-action is
approximately representable, also when $G$ is nonabelian.

Put $P=A^\alpha$ and $B=A\rtimes_\alpha G$, with canonical
unitaries $u_g\in B$. The expectations are
\[
 E_\alpha(a)=\frac1m\sum_g\alpha_g(a),\qquad
 F\!\left(\sum_ga_gu_g\right)=a_{1_G}.
\]
We explain the strict Rokhlin conclusion using the present results.
Pointwise outerness makes $B$ simple by
\cite[Theorem~3.1]{Kishimoto1981}. The nonzero Haar projection
$f=m^{-1}\sum_gu_g$ is full in $B$, so the action is saturated.
Proposition~\ref{prop:crossed-product-basic-construction} identifies
$P$ with the full corner $fBf$, and therefore $P$ is simple. It also
gives
\[
 C_r^*\langle A,e_P\rangle\cong B,\qquad
 F=\widehat E_\alpha,\qquad \operatorname{Index}(E_\alpha)=m.
\]
Theorem~\ref{thm:main} now implies that $E_\alpha$ is approximately
representable. Consequently $F=\widehat E_\alpha$ has the Rokhlin
property by \cite[Theorem~3.16(2)]{LO2019}. This is the strict Rokhlin
property, not just the weak tracial conclusion recorded in
\cite[Corollary~5.10]{LO2026WT}.

By \cite[Proposition~5.8]{LO2026WT}, the inclusion $A\subset B$,
with its specified expectation $F$, is not isomorphic to a fixed-point
inclusion arising from an ordinary finite-group action. The proof
there uses the noncocommutativity of $C(G)$, whereas every group
Hopf $\Cst$-algebra $C^*(L)$ is cocommutative. For the obstruction
argument, we refer to the cited proposition. This does not exclude a Hopf
fixed-point model: the dual action is
\[
 \delta_g\cdot\Bigl(\sum_ha_hu_h\Bigr)=a_gu_g,
 \qquad B^{C(G)}=A.
\]
Thus the example illustrates precisely the passage from ordinary
group actions to Hopf fixed-point inclusions.
\end{example}

The following example specializes the construction of
\cite[Section~7]{KT-2015}.
It shows that the action--expectation comparisons apply beyond
commutative and cocommutative Hopf $\Cst$-algebras.

\begin{example}\label{ex:kac-paljutkin}
Let $H=H_8$ be the eight-dimensional Kac--Paljutkin Hopf
$\Cst$-algebra \cite{KP1966}, which is neither commutative nor
cocommutative, and put $K=H_8^\circ$. Set
\[
 D=K\rtimes_{\Delta_K}H\cong M_8(\mathbb C),
 \qquad A=\bigotimes_{n=1}^\infty D.
\]
Here the action used to form $D$ is the one corresponding to the
coaction $\Delta_K\colon K\to K\otimes K$.
For dual bases $\{h_i\}$ and $\{h^i\}$, put
$V=\sum_i(1_K\rtimes h_i)\otimes h^i\in D\otimes K$ and
$u_n=V_{1,n+1}\cdots V_{n,n+1}$, where the last leg is $K$.
The compatible maps
\[
 \rho_n(a)=u_n(a\otimes1_K)u_n^*,\qquad a\in D^{\otimes n},
\]
induce the coaction $\rho\colon A\to A\otimes K$ of
\cite[Section~7]{KT-2015}. Let $H$ act on $A$ by slicing $\rho$.
By \cite[Propositions~7.4 and~7.5]{KT-2015}, this coaction is
approximately representable and has their Rokhlin property.
More explicitly, with $q=e_K\rtimes1_H\in D$, their tail projections
satisfy
\[
 p_n=1_D^{\otimes(n-1)}\otimes q,
 \qquad p=[(p_n)]\in A^\infty\cap A',
 \qquad E^\infty(p)=\tfrac18 1_A.
\]
Centrality follows because $p_n$ commutes with every fixed finite
tensor stage for all sufficiently large $n$.

We check the hypotheses of our comparisons. Applied to
$\lambda=\Sigma\rho$ in the convention fixed above,
\cite[Remark~6.2 and Theorem~4.6]{BSV2017} show that both $P=A^H$ and
$B=A\rtimes H$ are simple, since $A$ is a simple UHF algebra.
The nonzero Haar projection $f=j(e_H)\in B$ is therefore full, so
the action is saturated and $\operatorname{Index}(E)=8$.
The map $a\mapsto ap$ is a nonzero $*$-homomorphism on the simple
algebra $A$, and hence is injective. Thus $E$ has the Rokhlin
property, and Theorem~\ref{thm:Hopf-incl} gives our Rokhlin property
for the action. Theorem~\ref{thm:main} also gives approximate
representability of $E$. In particular, all the simplicity and
saturation hypotheses of Theorem~\ref{Rok-AR-duality} hold.
Since $K$ is noncommutative, this example cannot satisfy the older
full-central-range condition, by
\cite[Corollary~11.5]{KT-2015}.
\end{example}

The next consequences are already implied by the spatial approach:
\cite[Theorem~4.6]{BSV2017}, together with
\cite[Theorem~2.9]{BS2016}. We record them as applications of our
comparison, not as new permanence results. Here
$\dim_{\mathrm{nuc}}$ denotes nuclear dimension and
$\operatorname{tsr}$ denotes topological stable rank.

\begin{corollary}\label{cor:permanence}
Let $H\curvearrowright A$ be a saturated Rokhlin action of a
finite-dimensional Hopf $\Cst$-algebra on a separable unital simple
$\Cst$-algebra. Then
\[
 \max\{\dim_{\mathrm{nuc}}P,\dim_{\mathrm{nuc}}B\}
 \leq\dim_{\mathrm{nuc}}A.
\]
If $\operatorname{tsr}(A)=1$, then
$\operatorname{tsr}(P)=\operatorname{tsr}(B)=1$.
\end{corollary}

\begin{proof}
By Theorem~\ref{thm:RP-SRP} and
\cite[Theorem~4.6]{BSV2017}, the embeddings
\[
 P\longrightarrow A,
 \qquad B\longrightarrow M_N(A),\qquad N=\dim H,
\]
are sequentially split. The first splitting can also be read directly
from the fixed-point expectation: for a Rokhlin projection $p$,
the map $\beta(a)=N E^\infty(ap)\in P^\infty$ constructed in the
proof of Theorem~\ref{thm:RP-SRP} is a unital $*$-homomorphism and
$\beta(x)=x$ for $x\in P$.
Apply \cite[Theorem~2.9(3) and~(12)]{BS2016} to the two embeddings.
Nuclear dimension is unchanged by matrix amplification, and
$\operatorname{tsr}(A)=1$ implies $\operatorname{tsr}(M_N(A))=1$,
which gives the assertions.
\end{proof}

\subsection{Projection obstructions}

A unital $\Cst$-algebra is called \emph{projectionless} if $0$ and $1$ are its
only projections. The next result uses the projections that occur in the
two definitions for $E$. Saturation cannot be omitted in~(2): the trivial
action is approximately representable.

\begin{theorem}\label{thm:obstruction}
Let a finite-dimensional Hopf $\Cst$-algebra $H$ with $N = \dim H>1$ act on
a nonzero unital $\Cst$-algebra $A$.
\begin{enumerate}
 \item If $A$ is projectionless, then the action does not have the
       Rokhlin property.
 \item If the action is saturated and $P$ is projectionless, then
       neither the action nor $E\colon A\to P$ is approximately
       representable.
\end{enumerate}
Consequently, a saturated action of $H$ on a projectionless
$\Cst$-algebra has neither property.
\end{theorem}

\begin{proof}
For any $\Cst$-algebra $D$, a projection
$q\in D^\infty$ has bounded selfadjoint lifts $x_n$ with
$\|x_n^2-x_n\|\to0$. For all large $n$, functional calculus gives
projections $q_n=\chi_{[1/2,\infty)}(x_n)\in D$ with
$\|q_n-x_n\|\to0$; taking $q_n=0$ at the remaining indices gives
$q=[(q_n)]$.

(1) Let $\pi$ be a Rokhlin homomorphism and $p=\pi(e_K)$. The averaging
computation in the first part of the proof of Theorem~\ref{thm:Hopf-incl},
which uses only equivariance and unitality of $\pi$, gives
$E^\infty(p)=N^{-1}1_A$. Write $p=[(p_n)]$ with each $p_n$ a projection
in $A$. Since $p_n\in\{0,1_A\}$, we have $E(p_n)\in\{0,1_A\}$ and
therefore $\|E(p_n)-N^{-1}1_A\|\geq N^{-1}$ for every $n$, as
$N\geq2$. This is a contradiction.

(2) Suppose first that there are a projection $e\in P^\infty$ and
a finite set $\{u_i\}\subset A$ with
\begin{equation}\label{eq:obstruction}
 exe=E(x)e\quad(x\in A),\qquad \sum_iu_ieu_i^*=1.
\end{equation}
Write $e=[(e_n)]$ with each $e_n$ a projection in $P$, so that
$e_n\in\{0,1_A\}$. If $e_n=0$ for infinitely many $n$, then
$\|\sum_iu_ie_nu_i^*-1_A\|=1$ for those $n$, contradicting the second
relation. Hence $e_n=1_A$ for all large $n$, and $e=1$. The first relation
now reads $x=E(x)$ for every $x\in A$. Thus $A=P$ and
$E=\operatorname{id}_A$, so $\operatorname{Index}(E)=1_A$. But saturation
gives $\operatorname{Index}(E)=N1_A$ by
Proposition~\ref{prop:crossed-product-basic-construction}, a contradiction.

If $E$ is approximately representable, then \eqref{eq:obstruction} holds
by Definition~\ref{AR-incl}. If the action is approximately representable,
then \eqref{eq:obstruction} holds by Claims~2--4 in the proof of
Theorem~\ref{thm:main}; simplicity of $P$ is used there only for
condition~(3) of Definition~\ref{AR-incl}, which is not needed here.
This proves~(2). The last assertion follows because $P\subseteq A$ is a
unital subalgebra.
\end{proof}

\begin{example}\label{ex:jiang-su}
Let $\mathcal Z$ be the Jiang--Su algebra, which is simple, unital and
projectionless and satisfies $\mathcal Z\cong\bigotimes_{n=1}^\infty\mathcal Z$
\cite{JS1999}. Let $G$ be a nontrivial finite group. Put
$C=\bigotimes_{g\in G}\mathcal Z$, let $\beta_g$ permute its tensor factors
by left translation, and let
\[
 A=\bigotimes_{n=1}^\infty C\cong\mathcal Z,\qquad
 \alpha_g=\bigotimes_{n=1}^\infty\beta_g.
\]
This action is pointwise outer. Indeed, let $g\neq1_G$ and choose $c\in C$
with $\|c\|=1$ and $\delta=\|\beta_g(c)-c\|>0$. Let $c_n\in A$ be the copy of
$c$ in the $n$-th factor. If $\alpha_g=\operatorname{Ad}w$ with $w\in A$
unitary, then $\|wc_nw^*-c_n\|=\|wc_n-c_nw\|\to0$, because $(c_n)$ is a
central sequence. This contradicts $\|\alpha_g(c_n)-c_n\|=\delta$.

Put $H=C^*(G)$ and $K=C(G)$. By \cite[Theorem~3.1]{Kishimoto1981},
$B=A\rtimes_\alpha G$ is simple, so its Haar projection is full and
the action is saturated. Theorem~\ref{thm:obstruction} shows that
$\alpha$ has neither property. By Theorem~\ref{Rok-AR-duality}, the dual
action of $C(G)$ on $B$ has neither property. The latter conclusion is
not an instance of Theorem~\ref{thm:obstruction}(1), since $B$ contains the
nontrivial Haar projection. For nonabelian $G$ this is an action of a Hopf
$\Cst$-algebra that is not isomorphic to any $C^*(L)$ as a Hopf
$\Cst$-algebra.
\end{example}

\begingroup
\emergencystretch=1em
\printbibliography

@article{Izumi2004I,
  author       = {Izumi, Masaki},
  title        = {Finite group actions on {$C^*$}-algebras with the {R}ohlin
                  property. {I}},
  journal      = {Duke Math. J.},
  volume       = {122},
  number       = {2},
  pages        = {233--280},
  year         = {2004},
}

@article {JS1999,
    AUTHOR = {Jiang, Xinhui and Su, Hongbing},
     TITLE = {On a simple unital projectionless {$C^*$}-algebra},
   JOURNAL = {Amer. J. Math.},
    VOLUME = {121},
      YEAR = {1999},
    NUMBER = {2},
     PAGES = {359--413},
}

@article{Izumi2004II,
  author       = {Izumi, Masaki},
  title        = {Finite group actions on {$C^*$}-algebras with the {R}ohlin
                  property. {II}},
  journal      = {Adv. Math.},
  volume       = {184},
  number       = {1},
  pages        = {119--160},
  year         = {2004},
}

@incollection{OKT-2009,
  author       = {Osaka, Hiroyuki and Kodaka, Kazunori and Teruya, Tamotsu},
  title        = {The {R}ohlin property for inclusions of {$C^*$}-algebras
                  with a finite {W}atatani index},
  booktitle    = {Operator Structures and Dynamical Systems},
  series       = {Contemp. Math.},
  volume       = {503},
  pages        = {177--195},
  publisher    = {Amer. Math. Soc.},
  address      = {Providence, RI},
  year         = {2009},
  eprint       = {1001.4314},
  eprinttype   = {arXiv},
  eprintclass  = {math.OA},
}

@article{JP2008,
  author       = {Jeong, Ja A and Park, Gi Hyun},
  title        = {Saturated actions by finite dimensional {H}opf $*$-algebras
                  on {$C^*$}-algebras},
  journal      = {Internat. J. Math.},
  volume       = {19},
  number       = {2},
  pages        = {125--144},
  year         = {2008},
}

@article{SP1994,
  author       = {Peligrad, Costel and Szyma{\'n}ski, Wojciech},
  title        = {Saturated actions of finite dimensional {H}opf $*$-algebras
                  on {$C^*$}-algebras},
  journal      = {Math. Scand.},
  volume       = {75},
  pages        = {217--239},
  year         = {1994},
}

@article{KT-2015,
  author       = {Kodaka, Kazunori and Teruya, Tamotsu},
  title        = {The {R}ohlin property for coactions of finite dimensional
                  {$C^*$}-{H}opf algebras on unital {$C^*$}-algebras},
  year         = {2015},
  eprint       = {1209.4136},
  journal      = {Journal of Operator Theory},
  volume       = {74},
  number       = {2},
  pages        = {329--369},
  doi          = {10.7900/jot.2014jul02.2029},
}

@article{KT2009,
  author       = {Kodaka, Kazunori and Teruya, Tamotsu},
  title        = {Inclusions of unital {$C^*$}-algebras of index-finite type
                  with depth 2 induced by saturated actions of finite
                  dimensional {$C^*$}-{H}opf algebras},
  journal      = {Math. Scand.},
  volume       = {104},
  pages        = {221--248},
  year         = {2009},
}

@article{BSV2017,
  author       = {Barlak, Sel{\c c}uk and Szab{\'o}, G{\'a}bor and Voigt,
                  Christian},
  title        = {The spatial {R}okhlin property for actions of compact
                  quantum groups},
  journal      = {J. Funct. Anal.},
  volume       = {272},
  number       = {6},
  pages        = {2308--2360},
  year         = {2017},
  eprint       = {1605.08600},
  eprinttype   = {arXiv},
  eprintclass  = {math.OA},
}

@article{BS2016,
  author       = {Barlak, Sel{\c c}uk and Szab{\'o}, G{\'a}bor},
  title        = {Sequentially split $*$-homomorphisms between {$C^*$}-algebras},
  journal      = {International Journal of Mathematics},
  volume       = {27},
  number       = {13},
  pages        = {1650105},
  year         = {2016},
  doi          = {10.1142/S0129167X16501056},
  eprint       = {1510.04555},
  eprinttype   = {arXiv},
}

@article{KP1966,
  author       = {Kac, G. I. and Paljutkin, V. G.},
  title        = {Finite ring groups},
  journal      = {Trudy Moskovskogo Matematicheskogo Obshchestva},
  volume       = {15},
  pages        = {224--261},
  year         = {1966},
  langid       = {russian},
}

@misc{LO2026WT,
  author       = {Lee, Hyun Ho and Osaka, Hiroyuki},
  title        = {Rokhlin property and approximate representability
                  for inclusions of index-finite type},
  date         = {2026-08-19},
  eprint       = {1805.01613v3},
  eprinttype   = {arXiv},
  eprintclass  = {math.OA},
}

@article{Kishimoto1981,
  author       = {Kishimoto, Akitaka},
  title        = {Outer automorphisms and reduced crossed products
                  of simple {$C^*$}-algebras},
  journal      = {Communications in Mathematical Physics},
  volume       = {81},
  number       = {3},
  pages        = {429--435},
  year         = {1981},
  doi          = {10.1007/BF01209077},
}

@article{LO2019,
author = {Hyun Ho Lee and Hiroyuki Osaka},
title = {On dualities of actions and inclusions},
journal = {Journal of Functional Analysis},
volume = {276},
number = {2},
pages = {602-635},
year = {2019},
issn = {0022-1236},
doi = {https://doi.org/10.1016/j.jfa.2018.06.016},
url = {https://www.sciencedirect.com/science/article/pii/S0022123618302489},
}

@article{VD1997,
  author    = {Van Daele, Alfons},
  title     = {The Haar measure on finite quantum groups},
  journal   = {Proceedings of the American Mathematical Society},
  volume    = {125},
  number    = {12},
  pages     = {3489--3500},
  year      = {1997},
  doi       = {10.1090/S0002-9939-97-04037-9},
}

@book{timmermann,
    AUTHOR = {Timmermann, Thomas},
     TITLE = {An invitation to quantum groups and duality},
    SERIES = {EMS Textbooks in Mathematics},
      NOTE = {From Hopf algebras to multiplicative unitaries and beyond},
 PUBLISHER = {European Mathematical Society (EMS), Z\"urich},
      YEAR = {2008},
     PAGES = {xx+407},
      ISBN = {978-3-03719-043-2},
   MRCLASS = {46L05 (16W35 46L65 58B32)},
  MRNUMBER = {2397671},
MRREVIEWER = {Stefaan\ Vaes},
       DOI = {10.4171/043},
       URL = {https://doi.org/10.4171/043},
}

@article{BS1993,
  author  = {Baaj, Saad and Skandalis, Georges},
  title   = {Unitaires multiplicatifs et dualit{\'e} pour les produits crois{\'e}s de {$C^*$}-alg{\`e}bres},
  journal = {Ann. Sci. {\'E}cole Norm. Sup. (4)},
  volume  = {26},
  number  = {4},
  pages   = {425--488},
  year    = {1993},
  doi     = {10.24033/asens.1677},
}
\endgroup

\end{document}